\documentclass[a4paper,12pt]{amsart}

\usepackage{amsmath}
\usepackage{amssymb}
\usepackage{mathrsfs}
\usepackage{ifthen}
\usepackage{graphicx}
\usepackage[T1]{fontenc} 

\def\switchlinenumbers{\@ifstar
    {\let\makeLineNumberOdd\makeLineNumberRight
     \let\makeLineNumberEven\makeLineNumberLeft}%
    {\let\makeLineNumberOdd\makeLineNumberLeft
     \let\makeLineNumberEven\makeLineNumberRight}%
    }

\def\setmakelinenumbers#1{\@ifstar
  {\let\makeLineNumberRunning#1%
   \let\makeLineNumberOdd#1%
   \let\makeLineNumberEven#1}%
  {\ifx\c@linenumber\c@runninglinenumber
      \let\makeLineNumberRunning#1%
   \else
      \let\makeLineNumberOdd#1%
      \let\makeLineNumberEven#1%
   \fi}%
  }

\nonstopmode \numberwithin{equation}{section}
\newtheorem*{theorem*}{Theorem}

\newtheorem{thm}{Theorem}[section]
\newtheorem{cor}[equation]{Corollary}
\newtheorem{lem}[equation]{Lemma}
\newtheorem{prop}[equation]{Proposition}

\theoremstyle{definition}
\newtheorem{defn}{Definition}[section]

\newtheorem{qsn}[equation]{Question}
\newtheorem{prob}[equation]{Problem}
\newtheorem{rem}{Remark}[section]

\newtheorem{innercustomthm}{Theorem}
\newenvironment{customthm}[1]
  {\renewcommand\theinnercustomthm{#1}\innercustomthm}
  {\endinnercustomthm}

\newcounter{minutes}
\divide\time by 60
\newcounter{hours}
\multiply\time by 60
\addtocounter{minutes}{-\time}

\newcounter {own}
\def\theown {\thesection       .\arabic{own}}

\newenvironment{pf}[1][]{%
 \vskip 3mm
 \noindent
 \ifthenelse{\equal{#1}{}}%
  {{\slshape Proof. }}%
  {{\slshape #1.} }%
 }%
{\qed\bigskip}

\newcounter{alphabet}

\def\be{\begin{equation}}
\def\ee{\end{equation}}

\newcommand{\bee}{\begin{enumerate}}
\newcommand{\eee}{\end{enumerate}}

\newcommand{\blem}{\begin{lem}}
\newcommand{\elem}{\end{lem}}
\newcommand{\bthm}{\begin{thm}}
\newcommand{\ethm}{\end{thm}}
\newcommand{\bcor}{\begin{cor}}
\newcommand{\ecor}{\end{cor}}
\newcommand{\beg}{\begin{examp}}
\newcommand{\eeg}{\end{examp}}
\newcommand{\begs}{\begin{examples}}
\newcommand{\eegs}{\end{examples}}
\newcommand{\bdefe}{\begin{defin}}
\newcommand{\edefe}{\end{defin}}
\newcommand{\bprob}{\begin{prob}}
\newcommand{\eprob}{\end{prob}}
\newcommand{\bei}{\begin{itemize}}
\newcommand{\eei}{\end{itemize}}
\newcommand{\real}{{\operatorname{Re}\,}}

\newcommand{\norm}[1]{\left\lVert#1\right\rVert}
\newcommand{\abs}[1]{\left\lvert#1\right\rvert}
\newcommand{\innpdct}[1]{\left\langle#1\right\rangle}

\begin{document}

\title{Bohr radius for Banach spaces on simply connected domains}

\author{Vasudevarao Allu}
\address{Vasudevarao Allu,
School of Basic Sciences,
Indian Institute of Technology Bhubaneswar,
Bhubaneswar-752050, Odisha, India.}
\email{avrao@iitbbs.ac.in}

\author{Himadri Halder}
\address{Himadri Halder,
School of Basic Sciences,
Indian Institute of Technology Bhubaneswar,
Bhubaneswar-752050, Odisha, India.}
\email{himadrihalder119@gmail.com}

\subjclass[{AMS} Subject Classification:]{Primary 46E40, 47A56, 47A63; Secondary 46B20, 30B10, 30C20, 30C65}

\def\thefootnote{}
\footnotetext{ {\tiny File:~\jobname.tex,
printed: \number\year-\number\month-\number\day,
          \thehours.\ifnum\theminutes<10{0}\fi\theminutes }
} \makeatletter\def\thefootnote{\@arabic\c@footnote}\makeatother

\begin{abstract}
	Let $H^{\infty}(\Omega,X)$ be the space of bounded analytic functions $f(z)=\sum_{n=0}^{\infty} x_{n}z^{n}$ from a proper simply connected domain $\Omega$ containing the unit disk $\mathbb{D}:=\{z\in \mathbb{C}:|z|<1\}$ into a complex Banach space $X$ with $\norm{f}_{H^{\infty}(\Omega,X)} \leq 1$. Let $\phi=\{\phi_{n}(r)\}_{n=0}^{\infty}$ with $\phi_{0}(r)\leq 1$ such that $\sum_{n=0}^{\infty} \phi_{n}(r)$ converges locally uniformly with respect to $r \in [0,1)$. For $1\leq p,q<\infty$, we denote 
	\begin{equation*}
	R_{p,q,\phi}(f,\Omega,X)= \sup \left\{r \geq 0: \norm{x_{0}}^p \phi_{0}(r) + \left(\sum_{n=1}^{\infty} \norm{x_{n}}\phi_{n}(r)\right)^q \leq \phi_{0}(r)\right\}
	\end{equation*}
	and define the Bohr radius associated with $\phi$ by $$R_{p,q,\phi}(\Omega,X)=\inf \left\{R_{p,q,\phi}(f,\Omega,X): \norm{f}_{H^{\infty}(\Omega,X)} \leq 1\right\}.$$
	In this article, we extensively study the Bohr radius $R_{p,q,\phi}(\Omega,X)$, when $X$ is an arbitrary Banach space and $X=\mathcal{B}(\mathcal{H})$ is the algebra of all bounded linear operators on a complex Hilbert space $\mathcal{H}$. Furthermore, we establish the Bohr inequality for the operator-valued Ces\'{a}ro operator and Bernardi operator.
	
\end{abstract}

\maketitle
\pagestyle{myheadings}
\markboth{Vasudevarao Allu and  Himadri Halder}{Bohr radius for Banach spaces on simply connected domains}

\section{Introduction}
Let $H^{\infty}(\mathbb{D},\mathbb{C})$ be the space of bounded analytic functions from the unit disk $\mathbb{D}:=\{z \in \mathbb{C}:|z|<1\}$ into the complex plane $\mathbb{C}$ and we denote $\norm{f}_{\infty}:=\sup_{|z|<1} |f(z)|$. The remarkable theorem of Harald Bohr of a universal constant $r=1/3$ for functions in $H^{\infty}(\mathbb{D},\mathbb{C})$ is as follows.
\begin{customthm}{A}
	Let $f \in H^{\infty}(\mathbb{D},\mathbb{C})$ with the power series $f(z)=\sum_{n=0}^{\infty} a_{n}z^{n}$. If $\norm{f}_{\infty} \leq 1$, then 
	\begin{equation} \label{him-p7-e-1.1}
	\sum_{n=0}^{\infty} |a_{n}|r^{n} \leq 1
	\end{equation}
	for $|z|=r \leq 1/3$ and the constant $1/3$, referred to as the classical Bohr radius, is the best possible.
\end{customthm}
The Bohr's theorem has become popular when Dixon \cite{Dixon & BLMS & 1995} has used it to disprove a long-standing conjecture that if the non-unital von Neumann's inequality holds for a Banach algebra, then it is necessarily an operator algebra. 
It is important to note that \eqref{him-p7-e-1.1} can be written in the following equivalent form:
\begin{equation} \label{him-p7-e-1.2}
|a_{0}|\,\phi_{0}(r) + \sum_{n=1}^{\infty} |a_{n}|\, \phi_{n}(r) \leq \phi_{0}(r)
\end{equation}
for $r \leq R:= 1/3$, where $\phi_{n}(r)=r^n$ and $R$ is the smallest root of the equation $\phi_{0}(r)=2\sum_{n=1}^{\infty} \phi_{n}(r)$ in $(0,1)$. We observe that $\{\phi_{n}(r)\}_{n=0}^{\infty}$ is a sequence of non-negative continuous functions in $[0,1)$ such that the series $\sum_{n=0}^{\infty} \phi_{n}(r)$ converges locally uniformly with respect to $r \in [0,1)$. This fact leads to the following question. 
\begin{qsn}
Can we establish the inequality \eqref{him-p7-e-1.2} for any sequence $\{\psi_{n}(r)\}_{n=0}^{\infty}$ of non-negative continuous functions in $[0,1)$ such that the series $\sum_{n=0}^{\infty} \psi_{n}(r)$ converges locally uniformly with respect to $r \in [0,1)$.
\end{qsn}
We give the affirmative answer to this question in Theorem \ref{him-P7-thm-1.3}. 
In order to generalize the inequality \eqref{him-p7-e-1.2}, we first need to introduce some basic notations. Let $\mathcal{G}$ denote the set of all sequences $\phi= \{\phi_{n}(r)\}_{n=0}^{\infty}$ of non-negative continuous functions in $[0,1)$ such that the series $\sum_{n=0}^{\infty} \phi_{n}(r)$ converges locally uniformly with respect to $r \in [0,1)$. Now we want to define a modified Bohr radius associated with $\phi \in \mathcal{G}$.
\begin{defn} \label{him-vasu-P7-def-1.1}
Let $f \in H^{\infty}(\mathbb{D},\mathbb{C})$ with $f(z)=\sum_{n=0}^{\infty} a_{n}z^{n}$ such that $\norm{f}_{\infty} \leq 1$ in $\mathbb{D}$. For $\phi \in \mathcal{G}$, we denote 
\begin{equation} \label{him-p7-e-1.4}
R_{\phi}(f,\mathbb{C})= \sup \left\{r \geq 0: \sum_{n=0}^{\infty} |a_{n}|\phi_{n}(r) \leq \phi_{0}(r)\right\}.
\end{equation}
Define Bohr radius associated with $\phi$ by 
\begin{equation} \label{him-p7-e-1.5}
R_{\phi}(\mathbb{C})=\inf \left\{R_{\phi}(f,\mathbb{C}): \norm{f}_{\infty} \leq 1\right\}.
\end{equation}
\end{defn}
\noindent Clearly, $R_{\phi}(\mathbb{C})$ coincides with the classical Bohr radius $1/3$ for $\phi_{n}(r)=r^n$ for $r \in [0,1)$. In this article we are interested to study the operator-valued analogue of the Bohr radius $R_{\phi}(\mathbb{C})$, which we discuss in Definition \ref{him-vasu-P7-def-1.3}.
\vspace{3mm}

Over the past two decades there has been significant interest on several veriations of Bohr inequality \eqref{him-p7-e-1.1}. In $2000$, Djkaov and Ramanujan \cite{Djakov & Ramanujan & J. Anal & 2000} extensively studied the best possible constant $r_{p}$, for $1 \leq p < \infty $ such that 
\begin{equation} \label{him-p7-e-1.6}
\left(\sum_{n=0}^{\infty} |a_{n}|^p (r_{p})^{np}\right)^{1/p} \leq \norm{f}_{\infty},
\end{equation}
where $f(z)=\sum_{n=0}^{\infty} a_{n}z^n$. For $p=1$, $r_{p}$ coincides with the classical Bohr radius $1/3$. Using Haussdorf-Young's inequality, it is easy to see that $r_{p}=1$ for $p \in [2,\infty)$. Computing the precise value of $r_{p}$ for $1<p<2$ is difficult in general. This fact leads to estimate the value of $r_{p}$. The following best known estimate has been obtained in \cite{Djakov & Ramanujan & J. Anal & 2000}
\begin{equation} \label{him-p7-e-1.7}
\left(1+ \left(\frac{2}{p}\right)^{\frac{1}{2-p}}\right)^{\frac{p-2}{p}} \leq r_{p} \leq \inf \limits_{0 \leq a < 1} \frac{(1-a^p)^{1/p}}{\left((1-a^2)^p + a^p (1-a^p)\right)^{1/p}}.
\end{equation} 
 For further generalization of \eqref{him-p7-e-1.1} replacing $H^{\infty}$-norm by the $H^p$-norm, we refer to \cite{bene-2004}. 
Paulsen {\it et al.} \cite{paulsen-2002} have considered the another modification of \eqref{him-p7-e-1.1} and have shown that 
\begin{equation} \label{him-p7-e-1.9}
|a_{0}|^2 + \sum_{n=1}^{\infty} |a_{n}| \left(\frac{1}{2}\right)^n \leq 1,
\end{equation}
where $f(z)=\sum_{n=0}^{\infty} a_{n}z^n$ and $\norm{f}_{\infty} \leq 1$. Moreover, the constant $1/2$ is sharp. Using the same approach in \cite{paulsen-2002}, Blasco \cite{Blasco-OTAA-2010} has extended \eqref{him-p7-e-1.9} for the range of $p \in [1,2]$ and has shown that
\begin{equation} \label{him-p7-e-1.10}
|a_{0}|^p + \sum_{n=1}^{\infty} |a_{n}| \left(\frac{p}{p+2}\right)^n \leq 1.
\end{equation}
The constant $p/(p+2)$ is sharp. 
\vspace{3mm}

The study of Bohr radius has also been extended for functions defined on a proper simply connected domain of the complex plain. 
Throughout this paper, $\Omega$ stands for a simply connected domain containing the unit disk $\mathbb{D}$. Let $\mathcal{H}(\Omega)$ denote the class of analytic functions in $\Omega$ and let $\mathcal{B}(\Omega)$ be the class of functions $f \in \mathcal{H}(\Omega) $ such that $f(\Omega) \subseteq \overline{\mathbb{D}}$. The Bohr radius $B_{\Omega}$ for the class $\mathcal{B}(\Omega)$ is defined by (see \cite{Four-Rusc-2010})
$$
B_{\Omega}:=\sup\bigg\{r\in (0,1) : M_{f}(r)\leq 1\; \text{for all}\; f(z)=\sum_{n=0}^{\infty}a_nz^n\in\mathcal{B}(\Omega),\; z\in\mathbb{D}\bigg\},
$$
where $M_{f}(r):=\sum_{n=0}^{\infty}|a_n|r^n$ is the associated majorant series of $f \in \mathcal{B}(\Omega)$ in $\mathbb{D}$. It is easy to see that when $\Omega=\mathbb{D}$, $B_{\mathbb{D}}=1/3$, which is the classical Bohr radius for the class $\mathcal{B}(\mathbb{D})$.
\vspace{3mm}

For $0\leq \gamma<1$, we consider the following disk defined by 
$$
\Omega_{\gamma}:=\bigg\{z\in\mathbb{C} : \bigg|z+\frac{\gamma}{1-\gamma}\bigg|<\frac{1}{1-\gamma}\bigg\}.
$$
Clearly, $\Omega_{\gamma}$ contains $\mathbb{D}$ and $\Omega_{\gamma}$ reduces to $\mathbb{D}$ for $\gamma=0$. In $2010$,
Fournier and  Ruscheweyh  \cite{Four-Rusc-2010} studied Bohr inequality \eqref{him-p7-e-1.1} for the class $\mathcal{B}(\Omega_{\gamma})$.
\begin{thm}\cite{Four-Rusc-2010} \label{thm-1.2}
	For $ 0\leq \gamma<1 $, let $ f\in\mathcal{B}(\Omega_{\gamma}) $, with $ f(z)=\sum_{n=0}^{\infty}a_nz^n $ in $ \mathbb{D} $. Then,
	\begin{equation*}
	\sum_{n=0}^{\infty}|a_n|r^n\leq 1\;\; \text{for}\;\; r\leq\rho _{\gamma}:=\frac{1+\gamma}{3+\gamma}.
	\end{equation*}
	Moreover, $ \sum_{n=0}^{\infty}|a_n|\rho _{\gamma}^n=1 $ holds for a function $ f(z)=\sum_{n=0}^{\infty}a_nz^n $ in $ \mathcal{B}(\Omega_{\gamma}) $ if, and only if, $ f(z)=c $ with $ |c|=1 $.
\end{thm} 
\vspace{3mm}

The main aim of this paper is to study the vector-valued analogue of \eqref{him-p7-e-1.4}, \eqref{him-p7-e-1.5} and \eqref{him-p7-e-1.10} on simply connected domains and its connection with Banach space and Hilbert space theories. For discussing this, we first need to introduce some basic notations and give some definitions. Let $H^{\infty}(\mathbb{D},X)$ be the space of bounded analytic functions from $\mathbb{D}$ into a complex Banach space $X$ and we write $\norm{f}_{H^{\infty}(\mathbb{D},X)}= \sup_{|z|<1} \norm{f(z)}$. For $p \in [1, \infty)$, $H^{p}(\mathbb{D},X)$ denotes the space of analytic functions from $\mathbb{D}$ into $X$ such that 
\begin{equation} \label{him-p7-e-1.11}
\norm{f}_{H^{p}(\mathbb{D},X)}= \sup \limits _{0<r<1} \left(\int_{0}^{2\pi} \norm{f(re^{it})}^p \frac{dt}{2\pi}\right)^{1/p} < \infty.
\end{equation} 
Throughout this paper $\mathcal{B}(\mathcal{H})$ stands for the space of bounded linear operators on a complex Hilbert space $\mathcal{H}$. For any $T \in \mathcal{B}(\mathcal{H})$, $\norm{T}$ denotes the operator norm of $T$. Let $T \in \mathcal{B}(\mathcal{H})$. Then the adjoint operator $T^{*}:\mathcal{H} \rightarrow \mathcal{H}$ of $T$ defined by $\innpdct{Tx,y}=\innpdct{x,T^{*}y}$ for all $x, y \in \mathcal{H}$. $T$ is said to be normal if $T^{*}T=TT^{*}$, self-adjoint if $T^{*}=T$, and positive if $\innpdct{Tx,x} \geq 0$ for all $x \in \mathcal{H}$. The absolute value of $T$ is defined by $\abs{T}:=\left(T^{*}T\right)^{1/2}$, while $S^{1/2}$ denotes the unique positive square root of a positive operator $S$. Let $I$ be the identity operator on $\mathcal{H}$. 
\vspace{3mm}

Now we define the vector-valued analogue of Definition \ref{him-vasu-P7-def-1.1} on arbitrary simply connected domain containing the unit disk $\mathbb{D}$. Let $H^{\infty}(\Omega,X)$ be the space of bounded analytic functions from $\Omega$ into a complex Banach space $X$ and $\norm{f}_{H^{\infty}(\Omega,X)}= \sup_{z \in \Omega} \norm{f(z)}$.
\begin{defn} \label{him-vasu-P7-def-1.3}
	Let $f \in H^{\infty}(\Omega,X)$ be given by $f(z)=\sum_{n=0}^{\infty} x_{n}z^{n}$ in $\mathbb{D}$ with $\norm{f}_{H^{\infty}(\Omega,X)} \leq 1$. For $\phi \in \mathcal{G}$, we denote 
	\begin{equation} \label{him-p7-e-1.12}
	R_{\phi}(f,\Omega,X)= \sup \left\{r \geq 0: \sum_{n=0}^{\infty} \norm{x_{n}}\phi_{n}(r) \leq \phi_{0}(r)\right\}.
	\end{equation}
	Define Bohr radius associated with $\phi$ by 
	\begin{equation} \label{him-p7-e-1.13}
	R_{\phi}(\Omega,X)=\inf \left\{R_{\phi}(f,\Omega,X): \norm{f}_{H^{\infty}(\Omega,X)} \leq 1\right\}.
	\end{equation}
\end{defn}
It is important to note that for $\Omega=\Omega_{\gamma}$ and $\phi_{n}(r)=r^n$, by embedding $\mathbb{C}$ into $X$, from Theorem \ref{thm-1.2}, $R_{\phi}(\Omega_{\gamma},X) \leq (1+\gamma)/(3+\gamma)$ for every complex Banach space $X$. Clearly, $R_{\phi}(\mathbb{D},X)\leq 1/3$. However, this notion is not much significant in the finite-dimensional case for dimension greater than one. As usual, for $1 \leq p<\infty$, $\mathbb{C}^{m}_{p}$ stands for the space $\mathbb{C}^m$ endowed with the norm $\norm{w}_{p}=\left(\sum_{i=1}^{m} |w_{i}|^p\right)^{1/p}$ and $\norm{w}_{\infty}=\sup_{1\leq i \leq m} |w_{i}|$, where $w=(w_{1}, w_{2},\ldots,w_{m}) \in \mathbb{C}^m$. In \cite{Blasco-OTAA-2010}, Blasco has shown that $R_{\phi}(\mathbb{D},\mathbb{C}^{m}_{p})=0$ for $\phi_{n}=r^n$ in $[0,1)$ when $1\leq p\leq \infty$. By considering the same functions as in \cite{Blasco-OTAA-2010}, we show that, for $m \geq 2$, $R_{\phi}(\mathbb{D},\mathbb{C}^{m}_{p})$ need not be always non-zero for all $\phi \in \mathcal{G}$. In particular, we see that $R_{\phi}(\mathbb{D},\mathbb{C}^{m}_{p})$ becomes zero for some particular choices of $\phi$.
\begin{prop} \label{him-P7-prop-1.13}
	Let $\phi=\{\phi_{n}(r)\}^{\infty}_{n=0} \in \mathcal{G}$.
\begin{enumerate} 
	\item  For $m \geq 2$,  $R_{\phi}(\mathbb{D},\mathbb{C}^{m}_{\infty})=0$ when $r=0$ is the only zero of $\phi_{1}(r)$ in $[0,1)$.
 	\item For $1\leq p<\infty$ and $m \geq 2$,  $R_{\phi}(\mathbb{D},\mathbb{C}^{m}_{p})=0$ when $\phi_{0}(r)=1$ and $\phi_{1}(r)=\alpha \, r^{\beta}$ for $r \in [0,1)$ and $\alpha, \beta \in (0,\infty)$.
\end{enumerate}
\end{prop}

\begin{pf}
It is sufficient to prove for the case $m=2$.\\
(1) We consider the function $f(z)= (1,z)= e_{1} + e_{2}z$, $z \in \mathbb{D}$, where $e_{1}=(1,0)$ and $e_{2}=(0,1)$. Clearly, $\norm{f}_{H^{\infty}(\mathbb{D},\mathbb{C}^{2}_{\infty})}= \sup _{|z|<1}\norm{f(z)}_{\infty}=1$. Then from \eqref{him-p7-e-1.12}, we have 
\begin{equation*}
R_{\phi}(f,\mathbb{D},\mathbb{C}^{2}_{\infty}) = \sup \left\{r \geq 0: \norm{x_{0}}_{\infty} \phi_{0}(r) + \norm{x_{1}}_{\infty}\phi_{1}(r) \leq \phi_{0}(r)\right\},
\end{equation*}
where $x_{0}=e_{1}$ and $x_{1}=e_{2}$. Clearly, $\norm{x_{0}}_{\infty}=\norm{x_{1}}_{\infty}=1$. Then 
\begin{equation} \label{him-p7-e-1.16}
\norm{x_{0}}_{\infty} \phi_{0}(r) + \norm{x_{1}}_{\infty}\phi_{1}(r) = \phi_{0}(r) + \phi_{1}(r) \leq \phi_{0}(r),
\end{equation}
only when $\phi_{1}(r) \leq 0$ for $r \in [0,1)$. Thus, to obtain $R_{\phi}(f,\mathbb{D},\mathbb{C}^{2}_{\infty})$, we need to find the supremum of all such $r$ such that $\phi_{1}(r) \leq 0$ for $r \in [0,1)$. Since $\phi \in \mathcal{G}$, each $\phi_{n}(r)$ is non-negative for all $r \in [0,1)$. Therefore, \eqref{him-p7-e-1.16} holds only when $\phi_{1}(r)=0$ for $r \in [0,1)$. By the hypothesis, we have $\phi_{1}(r)=0$ if, and only if, $r=0$, which yields that \eqref{him-p7-e-1.16} holds only for $r=0$. Thus, $R_{\phi}(f,\mathbb{D},\mathbb{C}^{2}_{\infty})=0$ and so $R_{\phi}(\mathbb{D},\mathbb{C}^{2}_{\infty})=0$. This shows that $R_{\phi}(\mathbb{D},\mathbb{C}^{m}_{\infty})=0$.
\par (2) For $1<p<\infty$, using the fact $\lim _{s \rightarrow \infty} s^{1/p} - (s-1)^{1/p}=0$, for each $\epsilon>0$, one can easily find  a value $\delta \in (0,1)$ such that 
\begin{equation} \label{him-p7-e-1.17}
1-(1-\delta)^{1/p} < \alpha \, \epsilon^{\beta}\, \delta^{1/p}.
\end{equation}
We now consider the function 
\begin{equation*}
f(z)=\left((1-\delta)^{1/p}, \delta^{1/p}\,z\right)= (1-\delta)^{1/p} e_{1} + \delta^{1/p}\, e_{2} z.
\end{equation*}
It is easy to see that 
$$\norm{f} _{H^{\infty}(\mathbb{D},\mathbb{C}^{2}_{p})}=\sup_{|z|<1} \norm{f(z)}_{p}=\sup_{0<r<1} \left(((1-\delta)+ \delta r^{p})^{1/p}\right) =1,$$
and hence \eqref{him-p7-e-1.12} becomes 
\begin{equation} \label{him-p7-e-1.18}
R_{\phi}(f,\mathbb{D},\mathbb{C}^{2}_{p})= \sup \left\{r \geq 0: \norm{x_{0}}_{p} \phi_{0}(r) + \norm{x_{1}}_{p} \phi_{1}(r) \leq \phi_{0}(r)\right\}.
\end{equation}
In view of the assumptions $\phi_{0}(r)=1$ and $\phi_{1}(r)=\alpha r^{\beta}$, we have 
\begin{equation} \label{him-p7-e-1.19}
\norm{x_{0}}_{p} \phi_{0}(r) + \norm{x_{1}}_{p} \phi_{1}(r)= (1-\delta)^{1/p} + \delta^{1/p} \alpha r^{\beta}.
\end{equation}
Using \eqref{him-p7-e-1.19} in \eqref{him-p7-e-1.18}, we obtain
\begin{equation} \label{him-p7-e-1.20}
R_{\phi}(f,\mathbb{D},\mathbb{C}^{2}_{p})= \sup \{r \geq 0: (1-\delta)^{1/p} + \delta^{1/p} \alpha r^{\beta} \leq 1\}.
\end{equation}
Therefore, \eqref{him-p7-e-1.17} and \eqref{him-p7-e-1.20} show that $R_{\phi}(f,\mathbb{D},\mathbb{C}^{2}_{p}) \leq \epsilon$. Hence $R_{\phi}(\mathbb{D},\mathbb{C}^{2}_{p}) =0$ for $1<p<\infty$. Thus, $R_{\phi}(\mathbb{D},\mathbb{C}^{m}_{p}) =0$. \\[2mm]
Now for $p=1$, using the fact $\lim _{s \rightarrow \infty} \sqrt{s} - \sqrt{s-1}=0$, for each $\epsilon>0$, one can easily find  a value $\delta \in (0,1)$ such that 
\begin{equation} \label{him-p7-e-1.21}
1-\sqrt{1-\delta} < \alpha \, \epsilon^{\beta}\, \sqrt{\delta}.
\end{equation}
We consider the following function 
\begin{equation*}
f(z)=\frac{\sqrt{1-\delta}}{2} (1,1) + \frac{\sqrt{\delta}}{2}(1,-1)z = \frac{1}{2} \left(\sqrt{1-\delta} + \sqrt{\delta}z, \sqrt{1-\delta} - \sqrt{\delta}z\right).
\end{equation*}
A simple computation shows that 
\begin{align*}
\norm{f(z)}_{1} 
&= \frac{1}{2} \left(\abs{\sqrt{1-\delta} + \sqrt{\delta}z} + \abs{\sqrt{1-\delta} - \sqrt{\delta}z}\right)\\
& \leq \frac{1}{\sqrt{2}} \left(\abs{\sqrt{1-\delta} + \sqrt{\delta}z}^2 + \abs{\sqrt{1-\delta} - \sqrt{\delta}z}^2 \right)^{1/2}=1.
\end{align*}
By the similar lines of argument as above for the case $1<p<\infty$, we obtain $R_{\phi}(f,\mathbb{D},\mathbb{C}^{2}_{1}) \leq \epsilon$ and hence $R_{\phi}(\mathbb{D},\mathbb{C}^{2}_{1}) =0$. Thus, $R_{\phi}(\mathbb{D},\mathbb{C}^{m}_{1}) =0$.
\end{pf}
\begin{rem}
\begin{enumerate}
	\item If $\phi=\{\phi_{n}(r)\}_{n=0}^{\infty}$ with $\phi_{n}(r)=r^n$, then each $\phi_{n}$ is non-negative in $[0,1)$ and so $\phi \in \mathcal{G}$. Clearly, $\phi_{1}(r)=r$ has only zero at $r=0$ in $[0,1)$. In view of Proposition \ref{him-P7-prop-1.13} (1), the corresponding Bohr radius associated with $\phi$, is $R_{\phi}(\mathbb{D},\mathbb{C}^{m}_{\infty}) =0$. Furthermore, it is easy to see that $\phi_{0}(r)=1$ and $\phi_{1}(r)=\alpha\, r^{\beta}$ with $\alpha=\beta=1$ and hence by Proposition \ref{him-P7-prop-1.13} (2), we have $R_{\phi}(\mathbb{D},\mathbb{C}^{m}_{p}) =0$ for $1\leq p < \infty$ and $m \geq 2$.
	\item  Similarly, when $\phi=\{\phi_{n}(r)\}_{n=0}^{\infty}$ with $\phi_{n}(r)=(n+1)r^n, \, nr^n, \, n^2r^n$, Proposition \ref{him-P7-prop-1.13} gives the corresponding Bohr radius associated with $\phi$, $R_{\phi}(\mathbb{D},\mathbb{C}^{m}_{\infty}) =0$ and $R_{\phi}(\mathbb{D},\mathbb{C}^{m}_{p}) =0$ for $1\leq p < \infty$ and $m \geq 2$.
\end{enumerate}
\end{rem}
The above fact leads us to consider the vector-valued analogue of \eqref{him-p7-e-1.9} in a simply connected domain for a given Banach space $X$ and parameters $0<p,q<\infty$. We define a modified Bohr radius which need not be zero for all $\phi  \in \mathcal{G}$ even for infinite-dimensional Banach spaces.
\begin{defn} \label{him-vasu-P7-def-1.4}
	Let $f \in H^{\infty}(\Omega,X)$ be given by $f(z)=\sum_{n=0}^{\infty} x_{n}z^{n}$ in $\mathbb{D}$ with $\norm{f}_{H^{\infty}(\Omega,X)} \leq 1$. For $\phi=\{\phi_{n}(r)\}_{n=0}^{\infty}  \in \mathcal{G}$ with $\phi_{0}(r)\leq 1$, $1\leq p,q<\infty$, we denote 
	\begin{equation} \label{him-p7-e-1.22}
	R_{p,q,\phi}(f,\Omega,X)= \sup \left\{r \geq 0: \norm{x_{0}}^p \phi_{0}(r) + \left(\sum_{n=1}^{\infty} \norm{x_{n}}\phi_{n}(r)\right)^q \leq \phi_{0}(r)\right\}.
	\end{equation}
	Define Bohr radius associated with $\phi$ by 
	\begin{equation} \label{him-p7-e-1.23}
	R_{p,q,\phi}(\Omega,X)=\inf \left\{R_{p,q,\phi}(f,\Omega,X): \norm{f}_{H^{\infty}(\Omega,X)} \leq 1\right\}.
	\end{equation}
\end{defn}
Clearly, $R_{1,1,\phi}(\Omega,X)=R_{\phi}(\Omega,X)$. For $p_{1} \leq p_{2}$ and $q_{1} \leq q_{2}$, we have the following inclusion relation:
\begin{equation} \label{him-P7-e-1.22}
R_{p_{1},q_{1},\phi}(\Omega,X) \leq R_{p_{2},q_{2},\phi}(\Omega,X).
\end{equation}
Finding the exact value of $R_{p,q,\phi}(\Omega,X)$ is very difficult in general, even for $\Omega=\mathbb{D}$ and $X= \mathbb{C}^{1}_{2}$. In $2002$, Paulsen {\it et al.}\cite{paulsen-2002} proved that $R_{2,1,\phi}(\mathbb{D},\mathbb{C})=1/2$ for $\phi=\{\phi_{n}(r)\}^{\infty}_{n=0}$ with $\phi_{n}(r)=r^n$. Later, for the same $\phi$, Blasco \cite{Blasco-OTAA-2010} has shown that $R_{2,1,\phi}(\mathbb{D},\mathbb{C})=p/(p+2)$ for $1 \leq p \leq 2$. By considering the same example as in Proposition \ref{him-P7-prop-1.13}, we have the following interesting result.
\begin{prop} \label{him-P7-prop-1.23}
	Let $\phi=\{\phi_{n}(r)\}^{\infty}_{n=0} \in \mathcal{G}$.
	 For $m \geq 2$ and $1\leq p,q<\infty$,  $R_{p,q,\phi}(\mathbb{D},\mathbb{C}^{m}_{\infty})=0$ when $r=0$ is the only zero of $\phi_{1}(r)$ in $[0,1)$.
\end{prop}
It is important to note that $\mathbb{C}^m_{\infty}$ is not a Hilbert space. Indeed, let $x=(1,0,\ldots,0)$ and $y=(0,1,\ldots,0)$ be in $\mathbb{C}^m_{\infty}$. Then $\norm{x}_{\infty}=\norm{y}_{\infty}=\norm{x+y}_{\infty}= \norm{x-y}_{\infty}=1$ and $\norm{x+y}^2_{\infty} + \norm{x-y}^2_{\infty}=2 \neq 4=2\norm{x}^2_{\infty}+2\norm{y}^2_{\infty}$. Hence, Parallelogram law is violated.
 Blasco \cite{Blasco-OTAA-2010} has shown that for $m \geq 2$, $R_{p,p,\phi}(\mathbb{D},\mathbb{C}^m_{2})>0$ if, and only if, $p\geq 2$ when $\phi_{n}(r)=r^n$. It is worth mentioning that $X=\mathbb{C}^m_{2}$ is a Hilbert space with the inner product $\innpdct{.}$, where $\norm{x}_{2}=\sqrt{\innpdct{x,x}}$. This fact leads us to the following question.
  \begin{qsn} \label{him-P7-qsn-1.24}
  Does the radius $R_{p,q,\phi}(\mathbb{D},\mathcal{B}(\mathcal{H}))$ always "strictly" positive for $2\leq p \leq q$ under some suitable conditions on $\phi$ in case $X=\mathcal{B}(\mathcal{H})$, where $\mathcal{B}(\mathcal{H})$ is the space of bounded linear operators on a complex Hilbert space $\mathcal{H}$?
  \end{qsn}
We give an affirmative answer to the Question \ref{him-P7-qsn-1.24} in the following form.
\begin{thm} \label{him-P7-thm-1.2}
	Let $f \in H^{\infty}(\mathbb{D},\mathcal{B}(\mathcal{H}))$ be given by $f(z)= \sum_{n=0}^{\infty} A_{n}z^n$ in $\mathbb{D}$ with $A_{n} \in \mathcal{B}(\mathcal{H})$ for $n \in \mathbb{N} \cup \{0\}$ and $\norm{f(z)}_{H^{\infty}(\mathbb{D},\mathcal{B}(\mathcal{H}))} \leq 1$. Also let for $p\geq 2$, $\phi=\{\phi_{n}(r)\}_{n=0}^{\infty} \in \mathcal{G}$ satisfies the inequality
	\begin{equation} \label{him-P7-e-1.57}
	\phi_{0}(r)> 2 \sum_{n=1}^{\infty} \phi_{n}(r) \,\,\,\,\, \mbox{for} \,\,\,\, r\in [0,R(p)),  
	\end{equation}
	where $R(p)$ is the smallest root in $(0,1)$ of the equation
	\begin{equation} \label{him-P7-e-1.58}
	\phi_{0}(x)=2 \sum_{n=1}^{\infty} \phi_{n}(x).
	\end{equation}
	Then, for $p \geq 2$, we have $R_{p,p,\phi}(\mathbb{D},\mathcal{B}(\mathcal{H}))\geq R(p)$. That is, $R_{p,p,\phi}(\mathbb{D},\mathcal{B}(\mathcal{H}))>0$ for $p \geq 2$.
\end{thm}
\begin{pf}
	In view of the inclusion relation \eqref{him-P7-e-1.22}, it is enough to show that $R_{1,2,\phi}(\mathbb{D},\mathcal{B}(\mathcal{H}))>0$. By the given assumption, $f$ is in the unit ball of $H^{\infty}(\mathbb{D},\mathcal{B}(\mathcal{H}))$ {\it i.e.,}  $\norm{f}_{H^{\infty}(\mathbb{D},\mathcal{B}(\mathcal{H}))} \leq 1$. In particular, we have $\norm{f}^2_{H^{2}(\mathbb{D},\mathcal{B}(\mathcal{H}))} =\sum_{n=0}^{\infty} \norm{A_{n}}^2\leq 1$.
	Using Cauchy-Schwarz inequality, we obtain
	\begin{align} \label{him-p7-e-1.64}
	\norm{A_{0}} \phi_{0}(r) + \left(\sum_{n=1}^{\infty} \norm{A_{n}}\phi_{n}(r)\right)^2 
	& \leq \norm{A_{0}} \phi_{0}(r) + \left(\sum_{n=1}^{\infty} \norm{A_{n}}^2\right) \left(\sum_{n=1}^{\infty} \phi_{n}^2 (r)\right) \\ \nonumber
	& \leq \norm{A_{0}} \phi_{0}(r) + (1-\norm{A_{0}}^2) \sum_{n=1}^{\infty} \phi_{n}^2 (r) \leq \phi_{0}(r)\\ \nonumber
	& \leq \norm{A_{0}} \phi_{0}(r) + 2(1-\norm{A_{0}}) \sum_{n=1}^{\infty} \phi_{n}^2 (r) \leq \phi_{0}(r),
	\end{align}
	provided 
	\begin{equation} \label{him-p7-e-1.65}
	2 \sum_{n=1}^{\infty} \phi_{n}(r) < \phi_{0}(r).
	\end{equation}
	Now, by the given assumption \eqref{him-P7-e-1.57}, the inequality \eqref{him-p7-e-1.65} holds for $r \in [0,R(p))$, where $R(p)$ is the smallest root in $(0,1)$ of $\phi_{0}(r)=2 \sum_{n=1}^{\infty} \phi_{n}(r)$. Thus, we obtain that \eqref{him-p7-e-1.64} holds for $r \in [0,R(p))$. Hence $R_{1,2,\phi}(f,\mathbb{D},\mathcal{B}(\mathcal{H}))\geq R(p)$ and so, $R_{1,2,\phi}(\mathbb{D},\mathcal{B}(\mathcal{H}))\geq R(p)$. Since $R(p) \in (0,1)$, we have $R_{1,2,\phi}(\mathbb{D},\mathcal{B}(\mathcal{H}))>0$. Therefore, by the inclusion relation \eqref{him-P7-e-1.22}, for $p \geq 2$, we obtain $R_{p,p,\phi}(\mathbb{D},\mathcal{B}(\mathcal{H}))>0$. This completes the proof.
\end{pf}

\begin{rem}
By the virtue of the inclusion relation \eqref{him-P7-e-1.22} and Theorem \ref{him-P7-thm-1.2}, we conclude that $R_{p,q,\phi}(\mathbb{D},\mathcal{B}(\mathcal{H}))>0$ for $2\leq p \leq q$ under the same assumption on $\phi$ as in Theorem \ref{him-P7-thm-1.2}.
\end{rem}

As we have discussed the existence of the "strictly" positive radius $R_{p,q,\phi}(\mathbb{D},\mathcal{B}(\mathcal{H}))$ for $2\leq p \leq q$, it is natural to ask the following question.
\begin{qsn} \label{him-P7-qsn-1.34}
Does the radius $R_{p,q,\phi}(\mathbb{D},\mathcal{B}(\mathcal{H}))$ always strictly positive for $1\leq p,q <2$ under the some suitable conditions on $\phi$? If so, then is it possible to obtain the precise value of $R_{p,q,\phi}(\mathbb{D},\mathcal{B}(\mathcal{H}))$?
\end{qsn}
We give the affirmative answer to the Question \ref{him-P7-qsn-1.34}. We prove that $R_{p,q,\phi}(\mathbb{D},\mathcal{B}(\mathcal{H}))$ is strictly positive for $1\leq p,q <2$. Although finding the exact value of $R_{p,q,\phi}(\mathbb{D},\mathcal{B}(\mathcal{H}))$ for $1<q<2$ is very much complicated but we  we can find a good estimate of the Bohr radius $R_{p,1,\phi}(\Omega,\mathcal{B}(\mathcal{H}))$ on simply connected domain $\Omega$ containing $\mathbb{D}$. Let $f : \Omega \rightarrow \mathcal{B}(\mathcal{H})$ be a bounded analytic function {\it i.e.,} $f \in H^{\infty}(\Omega,\mathcal{B}(\mathcal{H}))$ with $f(z)=\sum_{n=0}^{\infty} A_{n}z^n$ in $\mathbb{D}$ such that $A_{n} \in \mathcal{B}(\mathcal{H})$ for all $n \in \mathbb{N} \cup \{0\}$. We denote
\begin{equation} \label{him-P7-e-1.24}
\lambda_{\mathcal{H}}:=\lambda_{\mathcal{H}}(\Omega):= \sup \limits _{\substack{f \in  H^{\infty}(\Omega,\mathcal{B}(\mathcal{H}))\\{\norm{f(z)}\leq1}}} \left\{\frac{\norm{A_{n}}}{\norm{\, I- |A_{0}|^2}\,} : A_{0} \not \equiv f(z)= \sum \limits_{n=0}^{\infty} A_{n} z^n, \,\, z \in \mathbb{D}\right\}.
\end{equation}
\begin{thm} \label{him-P7-thm-1.3}
For fixed $p \in [1,2]$. Let $f \in H^{\infty}(\Omega,\mathcal{B}(\mathcal{H}))$ be given by $f(z)= \sum_{n=0}^{\infty} A_{n}z^n$ in $\mathbb{D}$, where $A_{0}=\alpha_{0}I$ for $|\alpha_{0}|<1$ and $A_{n} \in \mathcal{B}(\mathcal{H})$ for all $n \in \mathbb{N} \cup \{0\}$ with $\norm{f}_{H^{\infty}(\Omega,\mathcal{B}(\mathcal{H}))} \leq 1$. If $\phi=\{\phi_{n}(r)\}_{n=0}^{\infty} \in \mathcal{G}$ satisfies the inequality
\begin{equation} \label{him-P7-e-1.25}
p \, \phi_{0}(r)> 2 \lambda_{\mathcal{H}} \sum_{n=1}^{\infty} \phi_{n}(r)\,\,\,\,\, \mbox{for} \,\,\,\, r \in [0, R_{\Omega}(p)),
\end{equation}
then the following inequality 
\begin{equation} \label{him-P7-e-1.26}
M_{f}(\phi, p, r):= \norm{A_{0}}^p \, \phi_{0}(r) + \sum_{n=1}^{\infty} \norm{A_{n}} \, \phi_{n}(r) \leq \phi_{0}(r)
\end{equation}
holds for $|z|=r \leq R_{\Omega}(p)$, where $R_{\Omega}(p)$ is the smallest root in $(0,1)$ of the equation 
\begin{equation} \label{him-P7-e-1.27}
p \, \phi_{0}(r)= 2 \lambda_{\mathcal{H}} \sum_{n=1}^{\infty} \phi_{n}(r).
\end{equation}
That is, $R_{\Omega}(p) \leq R_{p,1,\phi}(\Omega,\mathcal{B}(\mathcal{H}))$.
\end{thm}
\begin{pf}
Let $f \in H^{\infty}(\Omega,\mathcal{B}(\mathcal{H}))$ be given by $f(z)= \sum_{n=0}^{\infty} A_{n}z^n$ in $\mathbb{D}$ with $\norm{f(z)}_{H^{\infty}(\Omega,\mathcal{B}(\mathcal{H}))} \leq 1$. We note that $A_{0}=\alpha_{0}I$. Then, by the virtue of Lemma \ref{him-P7-e-1.24}, we have 
\begin{equation} \label{him-vasu-P6-e-2.59}
\norm{A_{n}} \leq \lambda_{\mathcal{H}}\norm{I- \abs{A_{0}^2}}= \lambda_{\mathcal{H}}\norm{I- \abs{\alpha_{0}}^2 I}=\lambda_{\mathcal{H}} (1- \abs{\alpha_{0}}^2)\, \,\,\,\, \mbox{for} \,\,\, n\geq 1.
\end{equation}
Using \eqref{him-vasu-P6-e-2.59}, we obtain 
\begin{align*}
M_{f}(\phi, p, r) 
& \leq \abs{\alpha_{0}}^p \, \phi_{0}(r) + \lambda_{\mathcal{H}} (1- \abs{\alpha_{0}}^2) \sum_{n=1}^{\infty} \phi_{n}(r) \\[2mm] 
&= \phi_{0}(r) + \lambda_{\mathcal{H}} (1- \abs{\alpha_{0}}^2) \left(\sum_{n=1}^{\infty} \phi_{n}(r) - \frac{(1-|\alpha_{0}|^p)}{\lambda_{\mathcal{H}}(1-|\alpha_{0}|^2)} \phi_{0}(r)\right).
\end{align*}
To obtain the inequality \eqref{him-P7-e-1.26} we now estimate the lower bound of $(1-|\alpha_{0}|^p)/\lambda_{\mathcal{H}} (1-|\alpha_{0}|^2)$. Let 
$$
B(x)=\frac{(1-x^p)}{\lambda_{\mathcal{H}}(1-x^2)} \,\,\,\,\, \mbox{for} \,\,\,\, x\in [0,1).
$$
For $p=2$, we have $B(x)=1/\lambda_{\mathcal{H}}$. For $p \in (0,2)$, let $\eta(x)=(2-p)x^p + px^{p-2} -2$. Then $B'(x)= -x \eta(x)/ (1-x^2)^2$ for $x \in (0,1)$. We note that $\eta'(x)= -p(2-p) x^{p-3} (1-x^2) <0$ for $x \in (0,1)$ and $p \in [1,2)$, which shows that $\eta$ is decreasing function in $(0,1)$ and thus, $\eta(x)>\eta(1)=0$ for $x \in (0,1)$. Therefore, $B'(x)<0$ in $(0,1)$ {\it i.e.,} $B$ is decreasing in $[0,1)$ and hence, 
$$
B(x) \geq \lim_{x \rightarrow 1^{-}} B(x)= \frac{p}{2\lambda_{\mathcal{H}}}\,\,\,\,\,\, \mbox{for} \,\,\,\,\, p\in [1,2).
$$
Thus, $B(x) \geq p/2\lambda_{\mathcal{H}}$ for $p \in [1,2]$, which leads to 
$$
M_{f}(\phi,p, r) \leq \phi_{0}(r) + \lambda_{\mathcal{H}} \left(1-\abs{\alpha_{0}}^2\right) \left(\sum_{n=1}^{\infty} \phi_{n}(r) - \frac{p}{2\lambda_{\mathcal{H}}} \phi_{0}(r)\right)
$$
and hence by \eqref{him-P7-e-1.25}, we obtain $M_{f}(\phi,p, r) \leq \phi_{0}(r)$ for $|z|=r \leq R_{\Omega}(p)$. Thus, $R_{\Omega}(p) \leq R_{p,1,\phi}(\Omega,\mathcal{B}(\mathcal{H}))$.
\end{pf}

When $p=1$ and $\phi_{n}(r)=r^n$, Theorem \ref{him-P7-thm-1.3} gives the following result, which is an analogue of classical Bohr inequality for operator valued analytic functions in a simply connected domain.
\begin{cor} \label{him-P7-cor-1.40}
Let $f \in H^{\infty}(\Omega,\mathcal{B}(\mathcal{H}))$ be given by $f(z)= \sum_{n=0}^{\infty} A_{n}z^n$ in $\mathbb{D}$, where $A_{0}=\alpha_{0}I$ for $|\alpha_{0}|<1$ and $A_{n} \in \mathcal{B}(\mathcal{H})$ for all $n \in \mathbb{N} \cup \{0\}$ with  $\norm{f(z)}_{H^{\infty}(\Omega,\mathcal{B}(\mathcal{H}))} \leq 1$. Then 
\begin{equation}
\sum_{n=0}^{\infty} \norm{A_{n}} r^n \leq 1 \,\,\,\,\,\, \mbox{for} \,\,\,\, r \leq \frac{1}{1+2\lambda_{\mathcal{H}}}.
\end{equation}
\end{cor}

As a consequence of Theorem \ref{him-P7-thm-1.3}, we wish to find the Bohr radius $R_{p,1,\phi}(\Omega_{\gamma},\mathcal{B}(\mathcal{H}))$ for the shifted disk $\Omega_{\gamma}$. For this, we need to compute the precise value of $\lambda_{\mathcal{H}}$, which in turns equivalent to study the coefficient estimates for the functions $f \in H^{\infty}(\Omega,\mathcal{B}(\mathcal{H}))$ of the form $f(z)= \sum_{n=0}^{\infty} A_{n}z^n$ in $\mathbb{D}$ with $\norm{f(z)}_{H^{\infty}(\Omega,\mathcal{B}(\mathcal{H}))} \leq 1$. To obtain the coefficient estimates, we shall make use of the following lemma from \cite{anderson-2006}.
\begin{lem} \cite{anderson-2006} \label{him-vasu-P6-lem-2.1}
	Let $B(z)$ be an analytic function with values in $\mathcal{B}(\mathcal{H})$ and satisfying $\norm{B(z)} \leq 1$ on $\mathbb{D}$. Then
	\begin{equation*}
	(1-|a|)^{n-1}\, \norm{\frac{B^{(n)}(a)}{n!}} \leq \frac{\norm{I-B(a)^{*}B(a)}^{1/2} \, \norm{I-B(a)B(a)^{*}}^{1/2}}{1-|a|^2} 
	\end{equation*}
	for each $a \in \mathbb{D}$ and $n=1,2,\ldots$.
\end{lem}
Using Lemma \ref{him-vasu-P6-lem-2.1}, we obtain the following coefficient estimates.
\begin{lem} \label{him-vasu-P6-lem-2.6}
	Let $f : \Omega_{\gamma}  \rightarrow \mathcal{B}(\mathcal{H})$ be analytic function with an expansion $f(z)= \sum_{n=0}^{\infty} A_{n}z^n$ in $\mathbb{D}$ such that $A_{n} \in \mathcal{B}(\mathcal{H})$ for all $n \in \mathbb{N} \cup \{0\}$ and $A_{0}$ is normal. Then 
	\begin{equation*}
	\norm{A_{n}} \leq \frac{\norm{I- |A_{0}|^2}}{1+\gamma} \,\,\,\,\, \mbox{for} \,\,\,n\geq 1.
	\end{equation*}
\end{lem}
\begin{pf}
	Let $\psi: \mathbb{D} \rightarrow \Omega_{\gamma}$ be analytic function defined by $\psi(z)= (z-\gamma)/(1-\gamma)$. Then, we see that the composition $g= f \circ \psi: \mathbb{D} \rightarrow \mathcal{B}(\mathcal{H})$ is analytic and 
	\begin{equation*}
	g(z)=f(\psi(z))= \sum \limits_{n=0}^{\infty} \frac{A_{n}}{(1-\gamma)^n}\, (z-\gamma)^n \,\,\,\,\,\,\, \mbox{for} \,\,\,\, |z-\gamma|<1-\gamma.
	\end{equation*}
	We note that $g(\gamma)=f(0)=A_{0}$ is normal and 
	\begin{equation} \label{him-vasu-P6-e-2.7}
	\frac{g^{(n)}(z)}{n!}= f^{(n)}\left(\frac{z-\gamma}{1-\gamma}\right) \,\, \frac{1}{(1-\gamma)^n}.
	\end{equation}
	In particular, for $z=\gamma$, \eqref{him-vasu-P6-e-2.7} gives 
	\begin{equation} \label{him-vasu-P6-e-2.8}
	(1-\gamma)^n \, \frac{g^{(n)}(\gamma)}{(n!)^2}= \frac{f^{(n)}(0)}{n!}=A_{n} \,\,\,\,\, \mbox{for} \,\, \, n \geq 1.
	\end{equation}
	In view of Lemma \ref{him-vasu-P6-lem-2.1}, we obtain 
	\begin{equation*}
	\norm{A_{n}} \leq (1-\gamma)^n \, \frac{g^{(n)}(\gamma)}{n!} \leq \frac{\norm{I-|g(\gamma)|^2}}{1+\gamma}=  \frac{\norm{I-|A_{0}|^2}}{1+\gamma} \,\,\,\,\, \mbox{for} \,\, \, n \geq 1.
	\end{equation*}
	This completes the proof.
\end{pf}

For $\Omega=\Omega_{\gamma}$, by making use of Lemma \ref{him-vasu-P6-lem-2.6} and \eqref{him-P7-e-1.24}, we obtain
\begin{equation}
\lambda_{\mathcal{H}}=\lambda_{\mathcal{H}}(\Omega_{\gamma})\leq\frac{1}{1+\gamma}.
\end{equation}
Now, we are in a position to find the Bohr radius $R_{p,1,\phi}(\Omega_{\gamma},\mathcal{B}(\mathcal{H}))$ for the shifted disk $\Omega_{\gamma}$.
\begin{thm} \label{him-vasu-P6-thm-2.8}
Fix $p \in [1,2]$. Let $f \in H^{\infty}(\Omega_{\gamma},\mathcal{B}(\mathcal{H}))$ be given by $f(z)= \sum_{n=0}^{\infty} A_{n}z^n$ in $\mathbb{D}$ with $\norm{f(z)}_{H^{\infty}(\Omega_{\gamma},\mathcal{B}(\mathcal{H}))} \leq 1$, where $A_{0}=\alpha_{0}I$ for $|\alpha_{0}|<1$ and $A_{n} \in \mathcal{B}(\mathcal{H})$ for all $n \in \mathbb{N} \cup \{0\}$. If $\phi=\{\phi_{n}(r)\}_{n=0}^{\infty} \in \mathcal{G}$ satisfies the inequality
\begin{equation} \label{him-vasu-P6-e-2.57}
\phi_{0}(r)> \frac{2}{p(1+\gamma)} \sum_{n=1}^{\infty} \phi_{n}(r) \,\,\,\,\, \mbox{for} \,\,\,\, r\in [0,R(p,\gamma)),  
\end{equation}
	then the inequality \eqref{him-P7-e-1.26} 
	holds for $|z|=r \leq R(p,\gamma)$, where $R(p,\gamma)$ is the smallest root in $(0,1)$ of the equation
	\begin{equation} \label{him-vasu-P6-e-2.56-a}
	\phi_{0}(x)=\frac{2}{p(1+\gamma)} \sum_{n=1}^{\infty} \phi_{n}(x).
	\end{equation}
	 Moreover, when $\phi_{0}(x)<(2/(p(1+\gamma))) \sum_{n=1}^{\infty} \phi_{n}(x)$ in some interval $\left(R(p,\gamma), R(p,\gamma)+ \epsilon\right)$ for $\epsilon >0$, then the constant $R(p,\gamma)$ cannot be improved further. That is, $R_{p,1,\phi}(\Omega_{\gamma},\mathcal{H}) = R(p,\gamma)$.
\end{thm}
\begin{pf}
For $\Omega=\Omega_{\gamma}$, $\lambda_{\mathcal{H}}=1/(1+\gamma)$, the condition \eqref{him-P7-e-1.25} becomes
\begin{equation*} 
\phi_{0}(r)> \frac{2}{p(1+\gamma)} \sum_{n=1}^{\infty} \phi_{n}(r) \,\,\,\,\, \mbox{for} \,\,\,\, r\in [0,R(p,\gamma)),  
\end{equation*}
where $R(p,\gamma)$ is the smallest root in $(0,1)$ of the equation
\begin{equation*}
\phi_{0}(x)=\frac{2}{p(1+\gamma)} \sum_{n=1}^{\infty} \phi_{n}(x).
\end{equation*}
By the virtue of Theorem \ref{him-P7-thm-1.3}, the required inequality \eqref{him-P7-e-1.26} holds for $r \in [0,R(p,\gamma))$. This gives that $R_{p,1,\phi}(\Omega_{\gamma},\mathcal{H})\geq R(p,\gamma)$. Our next aim is to show that $R_{p,1,\phi}(\Omega_{\gamma},\mathcal{H})= R(p,\gamma)$. For this, it is enough to show that the radius $R(p,\gamma)$ cannot be improved further. That is, $\norm{A_{0}}^p \, \phi_{0}(r) + \sum_{n=1}^{\infty} \norm{A_{n}} \, \phi_{n}(r) > \phi_{0}(r)$ holds for any $r>R(p,\gamma)$ {\it i.e.,} for any $r \in \left(R(p,\gamma), R(p,\gamma)+ \epsilon\right)$. To show this, we consider the following function
\begin{align} \label{him-P7-e-1.36-a}
F_a(z)=\left(\frac{a-\gamma-(1-\gamma)z}{1-a\gamma-a(1-\gamma)} \right)I\;\;\text{for}\;\; z\in \Omega_{\gamma}\;\; \text{and}\;\; a\in (0,1).
\end{align}
Define $ \psi_1 : \mathbb{D}\rightarrow\mathbb{D} $ by $ \psi_1(z)=(a-z)/(1-az) $ and $ \psi_2(z) :\Omega_{\gamma}\rightarrow\mathbb{D} $ by $ \psi_2(z)=(1-\gamma)z+\gamma $. 
Then, the function $ f_a=\psi_1\circ \psi_2 $ maps $ \Omega_{\gamma} $ univalently onto $ \mathbb{D}$. Thus, we note that $F_{a}(z)=f_{a}(z)I$ is analytic in $\Omega_{\gamma}$ and $\norm{F_{a}(z)} \leq \abs{f_{a}(z)} \leq 1$. A simple computation shows that
\begin{equation*}
F_a(z)=\left(\frac{a-\gamma-(1-\gamma)z}{1-a\gamma-a(1-\gamma)}\right)  I=A_0-\sum_{n=1}^{\infty}A_nz^n\;\; \mbox{for}\;\; z\in\mathbb{D},
\end{equation*} where $ a\in (0,1) $ and 
\begin{equation} \label{him-P7-e-1.37-a}
A_0=\frac{a-\gamma}{1-a\gamma}\, I\;\; \text{and}\;\; A_n=\left(\frac{1-a^2}{a(1-a\gamma)}\left(\frac{a(1-\gamma)}{1-a\gamma}\right)^n \right)I \,\,\,\,\, \mbox{for} \,\,\,\, n\geq 1.
\end{equation}
For the function $F_a$, we have 
\begin{align} \label{him-vasu-P6-e-2.60} &
\norm{A_{0}}^p \phi_{0}(r) + \sum_{n=1}^{\infty} \norm{A_{n}} \phi_{n}(r) \\[2mm] \nonumber
&= \left(\frac{a-\gamma}{1-a\gamma}\right)^p \phi_{0}(r) + (1-a^2) \sum_{n=1}^{\infty} \frac{a^{n-1} (1-\gamma)^n}{(1-a\gamma)^{n+1}} \phi_{n}(r) \\[2mm] \nonumber
&= \phi_{0}(r) + (1-a) \left(2\sum_{n=1}^{\infty} \phi_{n}(r) - p(1+\gamma) \phi_{0}(r)\right)+ \\[2mm] \nonumber
& \,\,\,\,\, (1-a) \left(\sum_{n=1}^{\infty}\frac{a^{n-1} (1+a) (1-\gamma)^n}{(1-a\gamma)^{n+1}} \phi_{n}(r) - 2 \sum_{n=1}^{\infty} \phi_{n}(r)\right)+  \\[2mm] \nonumber
& \,\,\,\, \left(p(1+\gamma)(1-a) + \left(\frac{a-\gamma}{1-a\gamma}\right)^p -1\right) \phi_{0}(r) \\[2mm] \nonumber
&= \phi_{0}(r) + (1-a) \left(2\sum_{n=1}^{\infty} \phi_{n}(r) - p(1+\gamma) \phi_{0}(r)\right) + O((1-a)^2)
\end{align}
as $a \rightarrow 1^{-}$. Also, we have that $2\sum_{n=1}^{\infty} \phi_{n}(r) > p(1+\gamma) \phi_{0}(r)$ for $r \in (R(p,\gamma), R(p,\gamma)+\epsilon)$. Then it is easy to see that the last expression of \eqref{him-vasu-P6-e-2.60} is strictly greater than $\phi_{0}(r)$ when $a$ is very close to $1$ {\it i.e.,} $a \rightarrow 1^{-}$ and $r \in (R(p,\gamma), R(p,\gamma)+\epsilon)$, which shows that the constant $R(p,\gamma)$ cannot be improved further. This completes the proof.
\end{pf}

The following are the consequences of Theorem \ref{him-vasu-P6-thm-2.8}.
\begin{cor}
	For $\psi_{n}(r)=r^n$ for $n \in \mathbb{N} \cup \{0\}$. Let $f$ be as in Theorem \ref{him-vasu-P6-thm-2.8}, then 
	\begin{equation} \label{him-vasu-P6-e-2.62}
	\norm{A_{0}}^p  + \sum_{n=1}^{\infty} \norm{A_{n}} r^n \leq 1 \,\,\,\,\, \mbox{for} \,\,\,\, |z|=r \leq R_{1}(p,\gamma):=\frac{p(1+\gamma)}{p(1+\gamma) +2}
	\end{equation}
	and the constant $R_{1}(p,\gamma)$ cannot be improved. Furthermore, if we consider complex valued analytic function $f \in \mathcal{B}(\Omega_{\gamma})$ such that $f(z)=\sum_{n=0}^{\infty} a_{n}z^n$ in $\mathbb{D}$ then from \eqref{him-vasu-P6-e-2.62}, we deduce that 
	\begin{equation} \label{him-vasu-P6-e-2.63}
	|a_{0}|^p  + \sum_{n=1}^{\infty} |a_{n}| r^n \leq 1 \,\,\,\,\, \mbox{for} \,\,\,\, |z|=r \leq R_{1}(p,\gamma):=\frac{p(1+\gamma)}{p(1+\gamma) +2}.
	\end{equation} 
\end{cor}
	We note that when $\Omega_{\gamma}=\mathbb{D}$ {\it i.e.,} $\gamma=0$, \eqref{him-vasu-P6-e-2.63} holds for $R_{1}(p):=p/(p+2)$, which has been independently obtained in \cite{Blasco-OTAA-2010}.
\begin{table}[ht]
	\centering
	\begin{tabular}{|l|l|l|l|l|}
		\hline
		$\gamma$& $R_{1}(1,\gamma)$& $R_{1}(1.5,\gamma)$& $R_{1}(1.7,\gamma)$& $R_{1}(2,\gamma)$ \\
		\hline
		$[0,0.2)$ & $[0.3333 \nearrow 0.3750)$& $[0.4285 \nearrow 0.4736)$ & $[0.4594 \nearrow 0.5050)$ & $[0.5000 \nearrow 0.5454)$\\
		\hline
		$[0.2,0.4)$& $[0.3750 \nearrow 0.4118)$& $[0.4736 \nearrow 0.5122)$& $[0.5050 \nearrow 0.5434)$& $[0.5454 \nearrow 0.5833)$\\
		\hline
		$[0.4,0.6)$& $[0.4118 \nearrow 0.4444)$& $[0.5122 \nearrow 0.5454)$& $[0.5434 \nearrow 0.5762)$& $[0.5833 \nearrow 0.6154)$\\
		\hline
		$[0.6,0.8)$& $[0.4444 \nearrow 0.4736)$& $[0.5454 \nearrow 0.5744)$& $[0.5762 \nearrow 0.6047)$& $[0.6154 \nearrow 0.6428)$\\
		\hline
		$[0.8,1)$& $[0.4736 \nearrow 0.5000)$& $[0.5744 \nearrow 0.6000)$& $[0.6047 \nearrow 0.6296)$& $[0.6428 \nearrow 0.6666)$\\
		
		\hline
	\end{tabular}
	\vspace{3mm}
	\caption{Values of $R_{1}(1,\gamma)$, $R_{1}(1.5,\gamma)$, $R_{1}(1.7,\gamma)$, and $R_{1}(2,\gamma)$ for various values of $\gamma \in [0,1)$.}
	\label{tabel-4.2-a}
\end{table}

\begin{figure}[!htb]
	\begin{center}
		\includegraphics[width=0.48\linewidth]{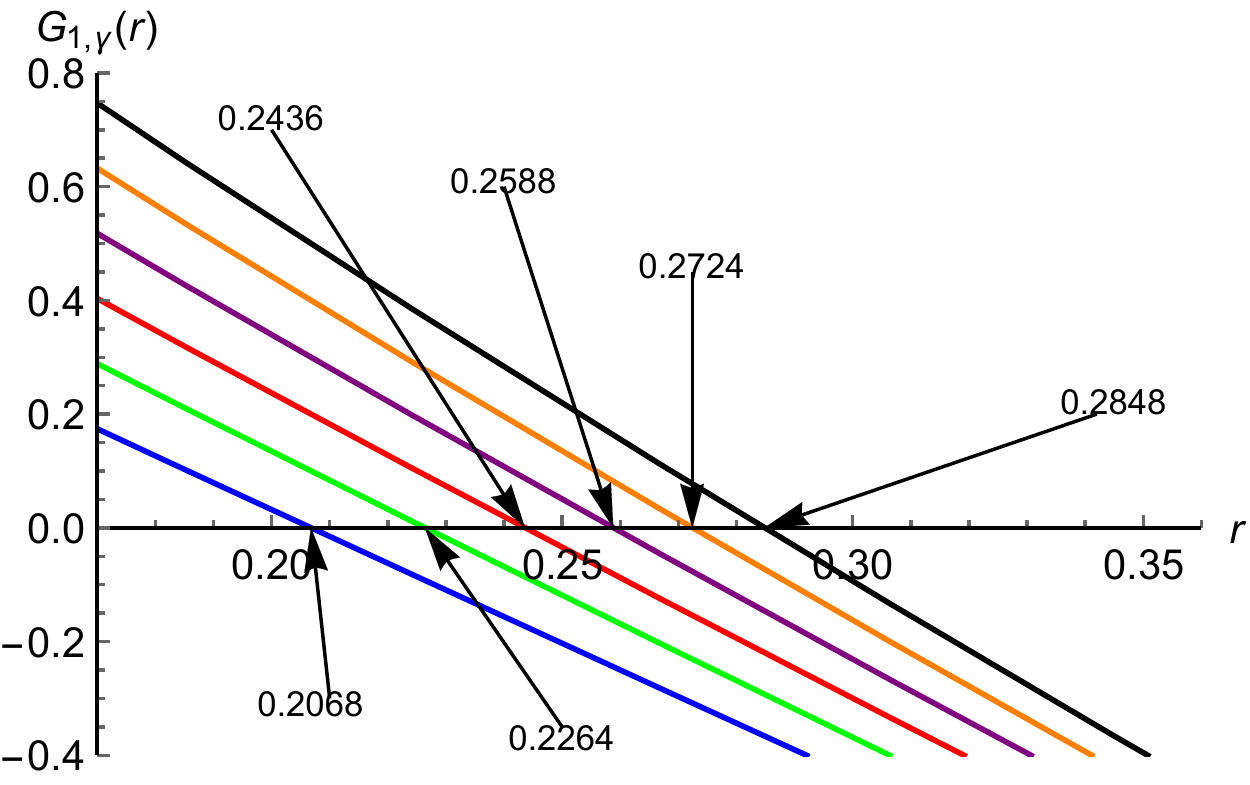}
		\,
		\includegraphics[width=0.48\linewidth]{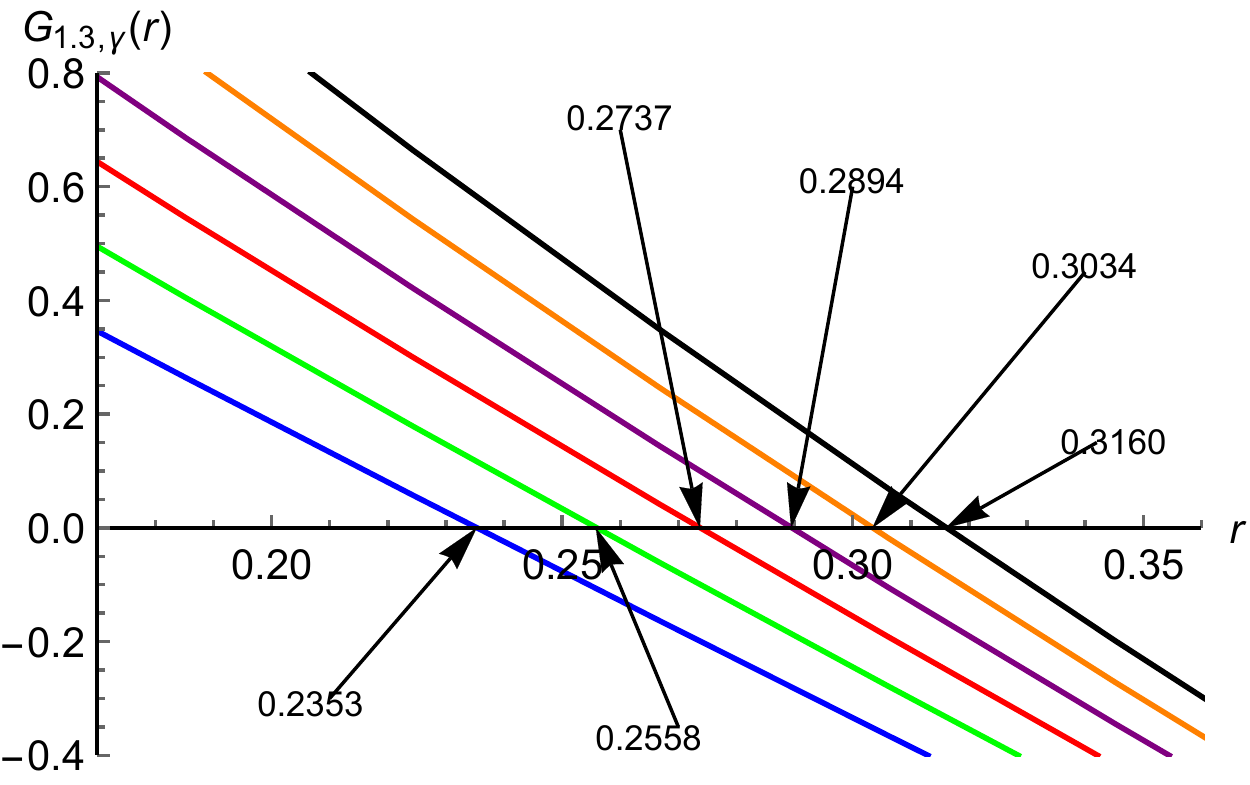}
	\end{center}
	\caption{The graph of $G_{1,\gamma}(r) $ and $G_{1.3,\gamma}( r)$ in $ (0,1) $ when $\gamma=0,0.2,0.4,0.6,0.8,1$.}
	\label{figure-4.2-a}
\end{figure}

\begin{figure}[!htb]
	\begin{center}
		\includegraphics[width=0.48\linewidth]{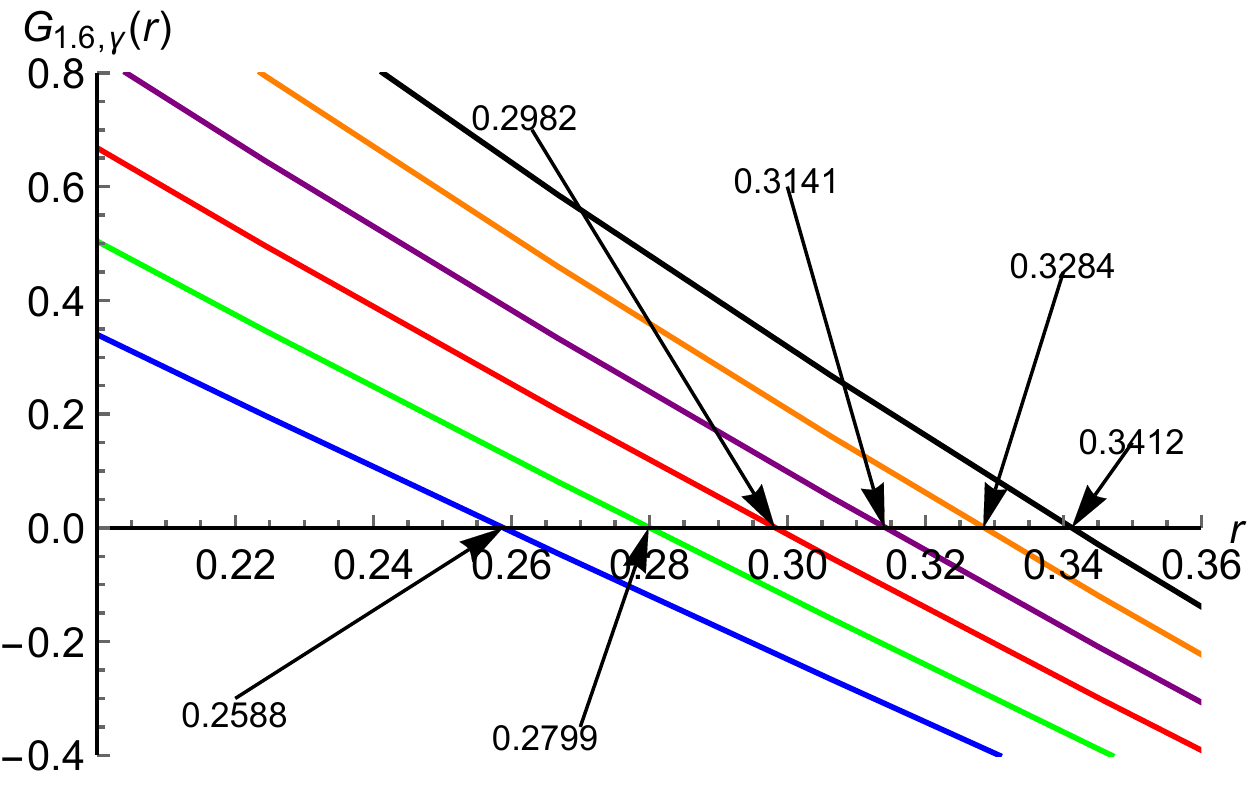}
		\,
		\includegraphics[width=0.48\linewidth]{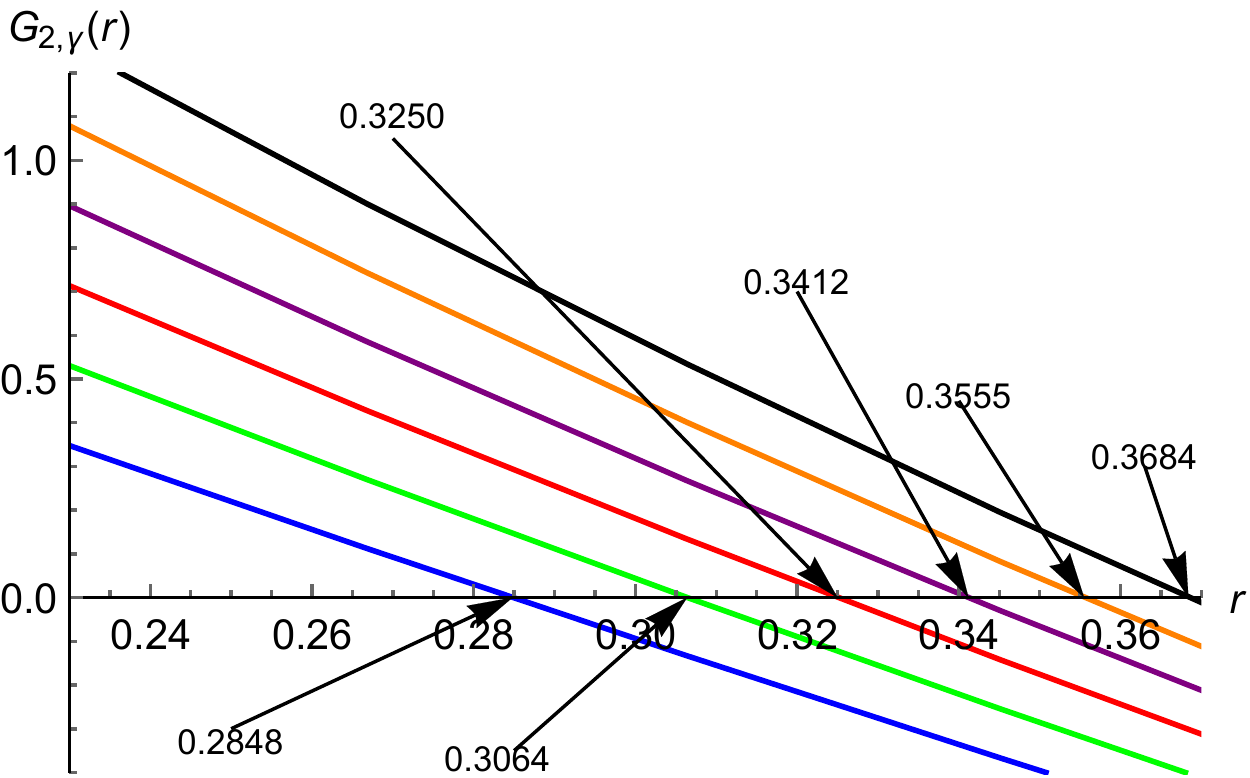}
	\end{center}
	\caption{The graph of $G_{1.6,\gamma}(r) $ and $G_{2,\gamma}( r)$ in $ (0,1) $ when $\gamma=0,0.2,0.4,0.6,0.8,1$.}
	\label{figure-4.2-b}
\end{figure}

\begin{cor}
	Let $\psi_{n}(r)=(n+1)r^n$ for $n \in \mathbb{N} \cup \{0\}$. Let $f$ be as in Theorem \ref{him-vasu-P6-thm-2.8}. Then we have the following sharp inequality 
	\begin{equation*} 
	\norm{A_{0}}^p  + \sum_{n=1}^{\infty} (n+1)\norm{A_{n}} r^n \leq 1 \,\,\,\,\, \mbox{for} \,\,\,\, |z|=r \leq R_{2}(p,\gamma):=1- \sqrt{\frac{2}{p(1+\gamma) +2}}.
	\end{equation*}
\end{cor}
\begin{table}[ht]
	\centering
	\begin{tabular}{|l|l|l|l|l|}
		\hline
		$\gamma$& $R_{1}(1,\gamma)$& $R_{1}(1.4,\gamma)$& $R_{1}(1.8,\gamma)$& $R_{1}(2,\gamma)$ \\
		\hline
		$[0,0.2)$ & $[0.1835 \nearrow 0.2094)$& $[0.2330 \nearrow 0.2628)$ & $[0.2745 \nearrow 0.3066)$ & $[0.2928 \nearrow 0.3258)$\\
		\hline
		$[0.2,0.4)$& $[0.2094 \nearrow 0.2330)$& $[0.2628 \nearrow 0.2893)$& $[0.3066 \nearrow 0.3348)$& $[0.3258 \nearrow 0.3545)$\\
		\hline
		$[0.4,0.6)$& $[0.2330 \nearrow 0.2546)$& $[0.2893 \nearrow 0.3132)$& $[0.3348 \nearrow 0.3598)$& $[0.3545 \nearrow 0.3798)$\\
		\hline
		$[0.6,0.8)$& $[0.2546 \nearrow 0.2745)$& $[0.3132 \nearrow 0.3348)$& $[0.3598 \nearrow 0.3822)$& $[0.3798 \nearrow 0.4023)$\\
		\hline
		$[0.8,1)$& $[0.2745 \nearrow 0.2928)$& $[0.3348 \nearrow 0.3545)$& $[0.3822 \nearrow 0.4023)$& $[0.4023 \nearrow 0.4226)$\\
		
		\hline
	\end{tabular}
	\vspace{3mm}
	\caption{Values of $R_{2}(1,\gamma)$, $R_{2}(1.4,\gamma)$, $R_{2}(1.8,\gamma)$, and $R_{2}(2,\gamma)$ for various values of $\gamma \in [0,1)$.}
	\label{tabel-4.2-b}
\end{table}
An observation shows that 
\begin{equation} \label{him-P7-e-1.55-a}
\sum_{n=1}^{\infty} nr^n= \frac{r}{(1-r)^2} \,\,\,\, \mbox{and} \,\,\, \sum_{n=1}^{\infty} n^2 r^n= \frac{r(1+r)}{(1-r)^3}. 
\end{equation}
Using \eqref{him-P7-e-1.55-a} and Theorem \ref{him-vasu-P6-thm-2.8}, we obtain the following corollary.
\begin{cor}
	Let $\psi_{0}(r)=1$ and $\psi_{n}(r)=n^k r^n$ for $n \geq 1$ and $k=1,2$.
	Then the following sharp inequalities hold 
	\begin{equation*}
	\norm{A_{0}}^p  + \sum_{n=1}^{\infty} n\norm{A_{n}} r^n \leq 1 \,\,\,\,\, \mbox{for} \,\,\,\, |z|=r \leq R_{3}(p,\gamma):= \frac{p(1+\gamma) + 1 - \sqrt{2p(1+\gamma)+ 1}}{p(1+\gamma) }
	\end{equation*}
	and 
	\begin{equation*}
	\norm{A_{0}}^p  + \sum_{n=1}^{\infty} n^2 \norm{A_{n}} r^n \leq 1 \,\,\,\,\, \mbox{for} \,\,\,\, |z|=r \leq R_{4}(p,\gamma),
	\end{equation*}
	where $R_{4}(p,\gamma)$ is the smallest positive root of the equation $G_{p,\gamma}(r):=p(1+\gamma)(1-r)^3 -2r(1+r)=0$ in $(0,1)$.
\end{cor}
\begin{table}[ht]
	\centering
	\begin{tabular}{|l|l|l|l|l|}
		\hline
		$\gamma$& $R_{3}(1,\gamma)$& $R_{3}(1.5,\gamma)$& $R_{3}(1.8,\gamma)$& $R_{1}(2,\gamma)$ \\
		\hline
		$[0,0.2)$ & $[0.2679 \nearrow 0.2967)$& $[0.3333 \nearrow 0.3640)$ & $[0.3640 \nearrow 0.3951)$ & $[0.3820 \nearrow 0.4132)$\\
		\hline
		$[0.2,0.4)$& $[0.2967 \nearrow 0.3218)$& $[0.3640 \nearrow 0.3903)$& $[0.3951 \nearrow 0.4216)$& $[0.4132 \nearrow 0.4396)$\\
		\hline
		$[0.4,0.6)$& $[0.3218 \nearrow 0.3441)$& $[0.3903 \nearrow 0.4132)$& $[0.4216 \nearrow 0.4444)$& $[0.4396 \nearrow 0.4624)$\\
		\hline
		$[0.6,0.8)$& $[0.3441 \nearrow 0.3640)$& $[0.4132 \nearrow 0.4334)$& $[0.4444 \nearrow 0.4645)$& $[0.4624 \nearrow 0.4823)$\\
		\hline
		$[0.8,1)$& $[0.3640 \nearrow 0.3820)$& $[0.4334 \nearrow 0.4514)$& $[0.4645 \nearrow 0.4823)$& $[0.4823 \nearrow 0.5000)$\\
		
		\hline
	\end{tabular}
	\vspace{3mm}
	\caption{Values of $R_{3}(1,\gamma)$, $R_{3}(1.5,\gamma)$, $R_{3}(1.8,\gamma)$, and $R_{3}(2,\gamma)$ for various values of $\gamma \in [0,1)$.}
	\label{tabel-4.2-c}
\end{table}
\begin{table}[ht]
	\centering
	\begin{tabular}{|l|l|l|l|l|}
		\hline
		$\gamma$& $R_{4}(1,\gamma)$& $R_{4}(1.3,\gamma)$& $R_{4}(1.6,\gamma)$& $R_{4}(2,\gamma)$ \\
		\hline
		$[0,0.2)$ & $[0.2068 \nearrow 0.2264)$& $[0.2353 \nearrow 0.2558)$ & $[0.2588 \nearrow 0.2799)$ & $[0.2848 \nearrow 0.3064)$\\
		\hline
		$[0.2,0.4)$& $[0.2264 \nearrow 0.2436)$& $[0.2558 \nearrow 0.2737)$& $[0.2799 \nearrow 0.2982)$& $[0.3064 \nearrow 0.3250)$\\
		\hline
		$[0.4,0.6)$& $[0.2436 \nearrow 0.2588)$& $[0.2737 \nearrow 0.2894)$& $[0.2982 \nearrow 0.3141)$& $[0.3250 \nearrow 0.3412)$\\
		\hline
		$[0.6,0.8)$& $[0.2588 \nearrow 0.2724)$& $[0.2894 \nearrow 0.3034)$& $[0.3141 \nearrow 0.3284)$& $[0.3412 \nearrow 0.3555)$\\
		\hline
		$[0.8,1)$& $[0.2724 \nearrow 0.2848)$& $[0.3034 \nearrow 0.3160)$& $[0.3284 \nearrow 0.3412)$& $[0.3555 \nearrow 0.3684)$\\
		
		\hline
	\end{tabular}
	\vspace{3mm}
	\caption{Values of $R_{4}(1,\gamma)$, $R_{4}(1.3,\gamma)$, $R_{4}(1.6,\gamma)$, and $R_{4}(2,\gamma)$ for various values of $\gamma \in [0,1)$.}
	\label{tabel-4.2-d}
\end{table}
From Table \ref{tabel-4.2-a}, Table \ref{tabel-4.2-b}, Table \ref{tabel-4.2-c}, and Table \ref{tabel-4.2-d}, for fixed values of $p$, we observe that Bohr radius $R_{1}(p,\gamma), R_{2}(p,\gamma), R_{3}(p,\gamma)$, and $R_{4}(p,\gamma)$ are monotonic increasing in $\gamma \in [0,1)$. In the Table \ref{tabel-4.2-a}, Table \ref{tabel-4.2-b}, Table \ref{tabel-4.2-c}, and Table \ref{tabel-4.2-d}, the notation $(R_{i}(p,\gamma_{1}) \nearrow R_{i}(p,\gamma_{2})]$ means that the value of $R_{i}(p,\gamma)$ is monotonically increasing from $\lim _{\gamma \rightarrow \gamma ^{+}_{1}}= R_{i}(\gamma_{1})$  to $ R_{i}(\gamma_{2})$ when $\gamma_{1} < \gamma \leq \gamma_{2}$, where $i=1,2,3$, and $4$. The Figure \ref{figure-4.2-a} and Figure \ref{figure-4.2-b} are devoted to the graphs of $G_{p,\gamma}(r)$ for different values of $p$ and $\gamma$. 

\section{Bohr inequality for Ces\'{a}ro operator}
In this section, we study the Bohr inequality for the operator-valued Ces\'{a}ro operator. For $\alpha \in \mathbb{C}$ with $\real \alpha >-1$, we have 
\begin{equation*} 
\frac{1}{(1-z)^{\alpha + 1}}= \sum_{k=0}^{\infty} C^{\alpha}_{k} z^k \,\,\,\,\,\, \mbox{where} \,\,\, C^{\alpha}_{k}= \frac{(\alpha +1) \ldots (\alpha + k)}{k!}.
\end{equation*}
Comparing the coefficient of $z^n$ on both the sides of the following identity
$$
\frac{1}{(1-z)^{\alpha + 1}} . \frac{1}{1-z}= \frac{1}{(1-z)^{\alpha + 2}},
$$
we obtain
\begin{equation} \label{him-vasu-P6-e-2.65}
C^{\alpha + 1}_{n} = \sum_{k=0}^{n} C^{\alpha}_{k} \,\,\, \, i.e.,\,\,\,\, \frac{1}{C^{\alpha + 1}_{n}} \sum_{k=0}^{n} C^{\alpha}_{k}=1.
\end{equation}
This property leads to consider the Ces\'{a}ro operator of order $\alpha$ or $\alpha$-Ces\'{a}ro operator (see \cite{stempak-1994}) on the space $\mathcal{H}(\mathbb{D})$ of analytic functions $f(z)= \sum_{n=0}^{\infty} a_{n} z^n$ in $\mathbb{D}$, which is defined by 
\begin{equation} \label{him-vasu-P6-e-2.66}
\mathcal{C}^{\alpha} f (z):= \sum_{n=0}^{\infty} \left(\frac{1}{C^{\alpha + 1}_{n}} \sum_{k=0}^{n} C^{\alpha}_{k}\, a_{k}\right) z^n. 
\end{equation}
A simple computation with power series gives the following integral form (see \cite{stempak-1994})
\begin{equation} \label{him-vasu-P6-e-2.67}
\mathcal{C}^{\alpha} f (z):= (\alpha + 1)\int\limits_{0}^{1} f(tz) \, \frac{(1-t)^{\alpha}}{(1-tz)^{\alpha +1}}\, dt
\end{equation}
with $\real \alpha > -1$. For $\alpha=0$, \eqref{him-vasu-P6-e-2.66} and \eqref{him-vasu-P6-e-2.67} give the classical Ces\'{a}ro operator 
\begin{equation} \label{him-vasu-P6-e-2.68}
\mathcal{C}f(z):= \mathcal{C}^{0} f (z)=\sum_{n=0}^{\infty} \left(\frac{1}{n+1} \sum_{k=0}^{n} a_{k}\right)z^n= \int\limits_{0}^{1} \frac{f(tz)}{1-tz} \, dt, \,\,\,\,\,\, z \in \mathbb{D}.
\end{equation}
In 1932, Hardy and Littlehood \cite{Hardy-1931} considered the classical Ces\'{a}ro operator and later, several authors have studied the boundedness of this operator on various function spaces (see \cite{albanese-2018}). In 2020, Berm\'{u}dez {\it et al.} \cite{Bermudez-2020} extensively studied the Ces\'{a}ro mean and boundedness of Ces\'{a}ro operators on Banach spaces and Hilbert spaces. 
\vspace{3mm}

In the same spirit of the definitions \eqref{him-vasu-P6-e-2.66} and \eqref{him-vasu-P6-e-2.67}, we define the Ces\'{a}ro operator on the space of analytic functions $f : \mathbb{D} \rightarrow \mathcal{B}(\mathcal{H})$ by 
\begin{equation} \label{him-vasu-P6-e-2.69}
\mathcal{C}^{\alpha}_{\mathcal{H}} f (z):= \sum_{n=0}^{\infty} \left(\frac{1}{C^{\alpha + 1}_{n}} \sum_{k=0}^{n} C^{\alpha}_{k}\, A_{k}\right) z^n = (\alpha + 1)\int\limits_{0}^{1} f(tz) \, \frac{(1-t)^{\alpha}}{(1-tz)^{\alpha +1}}\, dt, 
\end{equation}
where $f(z)=\sum_{n=0}^{\infty} A_{n}z^n$ in $\mathbb{D}$ and $A_{n},\, B_{n} \in \mathcal{B}(\mathcal{H})$ for all $n \in \mathbb{N} \cup \{0\}$. In \cite{Kayumov-Ponnuswamy-Cesaro-CRM-2020} and \cite{Kayumov-Ponnusamy-Cesaro-average-arxiv-2021}, Kayumov {\it et al.} have established an analogue of the Bohr theorem for the classical Ces\'{a}ro operator $\mathcal{C}f(z)$ and $\alpha$-Ces\'{a}ro operator $\mathcal{C}^{\alpha}_{\mathcal{H}}f(z)$ respectively. For an analytic function $f : \mathbb{D} \rightarrow \mathcal{B}(\mathcal{H})$ with $f(z)=\sum_{n=0}^{\infty} A_{n}z^n$ in $\mathbb{D}$, where $A_{n},\, B_{n} \in \mathcal{B}(\mathcal{H})$ for all $n \in \mathbb{N} \cup \{0\}$, we define the Bohr's sum by 
\begin{equation} \label{him-vasu-P6-e-2.70}
\mathcal{C}^{\alpha}_{f}(r):= \sum_{n=0}^{\infty} \left(\frac{1}{C^{\alpha + 1}_{n}} \sum_{k=0}^{n} C^{\alpha}_{k}\, \norm{A_{k}}\right)r^n \,\,\,\,\,\, \mbox{for} \,\,\, |z|=r.
\end{equation} 
Now we establish the counterpart of the Bohr theorem for $\mathcal{C}^{\alpha}_{\mathcal{H}} f (z)$.

\begin{thm} \label{him-vasu-P6-thm-2.9}
	Let $f : \Omega_{\gamma} \rightarrow \mathcal{B}(\mathcal{H})$ be an analytic function with $\norm{f(z)} \leq 1$ in $\Omega_{\gamma}$ such that $f(z)= \sum_{n=0}^{\infty} A_{n}z^n$ in $\mathbb{D}$, where $A_{0}=\alpha_{0}I$ for $|\alpha_{0}|<1$ and $A_{n} \in \mathcal{B}(\mathcal{H})$ for all $n \in \mathbb{N} \cup \{0\}$. Then for $\alpha>-1$, we have
	\begin{equation} \label{him-vasu-P6-e-2.70-a}
	\mathcal{C}^{\alpha}_{f}(r) \leq (\alpha + 1) \sum_{n=0}^{\infty} \frac{r^n}{\alpha + n +1}= \frac{(\alpha + 1)}{r^{\alpha +1}} \int\limits_{0}^{r} \frac{t^{\alpha}}{1-t} \, dt
	\end{equation}
	for $|z|=r \leq R(\gamma,\alpha)$, where $R(\gamma,\alpha)$ is the smallest root in $(0,1)$ of $C_{\gamma,\alpha}(r)=0$, where 
	\begin{equation*}
	C_{\gamma,\alpha}(r)=(3+\gamma)(1+\alpha)\sum_{n=0}^{\infty} \frac{r^n}{\alpha + n +1} -\frac{2}{1-r}.
	\end{equation*}
	The constant $R(\gamma,\alpha)$ cannot be improved further. 
\end{thm}

\begin{pf}
	Let $\alpha$-Ces\'{a}ro operator $\mathcal{C}^{\alpha}_{\mathcal{H}}f(z)$ be expressed in the following equivalent form 
	\begin{equation} \label{him-vasu-P6-e-2.71}
	\mathcal{C}^{\alpha}_{\mathcal{H}}f(z)= \sum_{n=0}^{\infty} A_{n} \phi_{n}(z),
	\end{equation}
	where $\phi_{n}(z)$ can be obtained by collecting the terms involving only $A_{n}$ in the right hand side of \eqref{him-vasu-P6-e-2.69}. Then it is easy to see that
	\begin{equation} \label{him-vasu-P6-e-2.73}
	\phi_{n}(z)= \sum_{k=n}^{\infty} \frac{C^{\alpha}_{k-n}}{C^{\alpha +1}_{k}} z^k
	\end{equation} 
	and hence by using the definition of $C^{\alpha}_{k}$, for $\alpha$-Ces\'{a}ro operator $\mathcal{C}^{\alpha}_{\mathcal{H}}f(z)$, we obtain
	\begin{equation} \label{him-vasu-P6-e-2.74}
	\phi_{0}(z)= \sum_{k=0}^{\infty} \frac{C^{\alpha}_{k}}{C^{\alpha +1}_{k}} z^k = (\alpha + 1) \sum_{k=0}^{\infty} \frac{z^k}{k+\alpha+1}, \,\,\,\,\, z \in \mathbb{D}.
	\end{equation} 
	It is easy to see that 
	\begin{equation} \label{him-vasu-P6-e-2.75}
	\mathcal{C}^{\alpha}_{f}(r)= \sum_{n=0}^{\infty} \norm{A_{n}} \phi_{n}(r).
	\end{equation}
By setting $f(z)=f_{1}(z):=(1/(1-z))I$ in \eqref{him-vasu-P6-e-2.71}, using \eqref{him-vasu-P6-e-2.65} and \eqref{him-vasu-P6-e-2.69}, we obtain 
	\begin{equation} \label{him-vasu-P6-e-2.76}
	\sum_{n=0}^{\infty} I \phi_{n}(z) = \mathcal{C}^{\alpha}_{\mathcal{H}}f_{1}(z)= \left(\frac{1}{1-z}\right)I.
	\end{equation}
	By using \eqref{him-vasu-P6-e-2.75} and \eqref{him-vasu-P6-e-2.76}, we obtain 
	\begin{equation} \label{him-vasu-P6-e-2.77}
	\mathcal{C}^{\alpha}_{f}(r)=\sum_{n=0}^{\infty} I \phi_{n}(r) = \frac{1}{1-r}.
	\end{equation}
	Thus, \eqref{him-vasu-P6-e-2.56-a} with $p=1$ takes the following form 
	\begin{equation*}
	(\alpha+1)\sum_{k=0}^{\infty} \frac{x^k}{k+\alpha+1}= \frac{2}{1+\gamma} \left(\frac{1}{1-x} - (\alpha+1)\sum_{k=0}^{\infty} \frac{x^k}{k+\alpha+1}\right),
	\end{equation*}
	which is equivalently,
	\begin{equation*}
	(3+\gamma)(\alpha+1)\sum_{k=0}^{\infty} \frac{x^k}{k+\alpha+1}= \frac{2}{1-x}.
	\end{equation*}
	Now the inequality \eqref{him-vasu-P6-e-2.70-a} follows from Theorem \ref{him-vasu-P6-thm-2.8}. Sharpness part follows from Theorem \ref{him-vasu-P6-thm-2.8}. This completes the proof.
\end{pf}
\begin{table}[ht]
	\centering
	\begin{tabular}{|l|l|l|l|l|}
		\hline
		$\gamma$& $R(\gamma,0)$& $R(\gamma,10)$& $R(\gamma,20)$& $R(\gamma,30)$ \\
		\hline
		$[0,0.3)$ & $[0.5335 \nearrow 0.6054)$& $[0.9860 \nearrow 0.9876)$ & $[0.9937 \nearrow 0.9943)$ & $[0.9961 \nearrow 0.9966)$\\
		\hline
		$[0.3,0.5)$& $[0.6054 \nearrow 0.6434)$& $[0.9876 \nearrow 0.9885)$& $[0.9945 \nearrow 0.9949)$& $[0.9966 \nearrow 0.9968)$\\
		\hline
		$[0.5,0.7)$& $[0.6434 \nearrow 0.6756)$& $[0.9885 \nearrow 0.9892)$& $[0.9949 \nearrow 0.9952)$& $[0.9968 \nearrow 0.9970)$\\
		\hline
		$[0.7,0.9)$& $[0.6756 \nearrow 0.7031)$& $[0.9892 \nearrow 0.9899)$& $[0.9952 \nearrow 0.9955)$& $[0.9970 \nearrow 0.9972)$\\
		\hline
		$[0.9,1)$& $[0.7031 \nearrow 0.7153)$& $[0.9899 \nearrow 0.9902)$& $[0.9955 \nearrow 0.9956)$& $[0.9972 \nearrow 0.9973)$\\
		
		\hline
	\end{tabular}
	\vspace{3mm}
	\caption{Values of $R(\gamma,0)$, $R(\gamma,10)$, $R(\gamma,20)$, and $R(\gamma,30)$ for various values of $\gamma \in [0,1)$.}
	\label{tabel-4.2-e}
\end{table}
\begin{figure}[!htb]
	\begin{center}
		\includegraphics[width=0.48\linewidth]{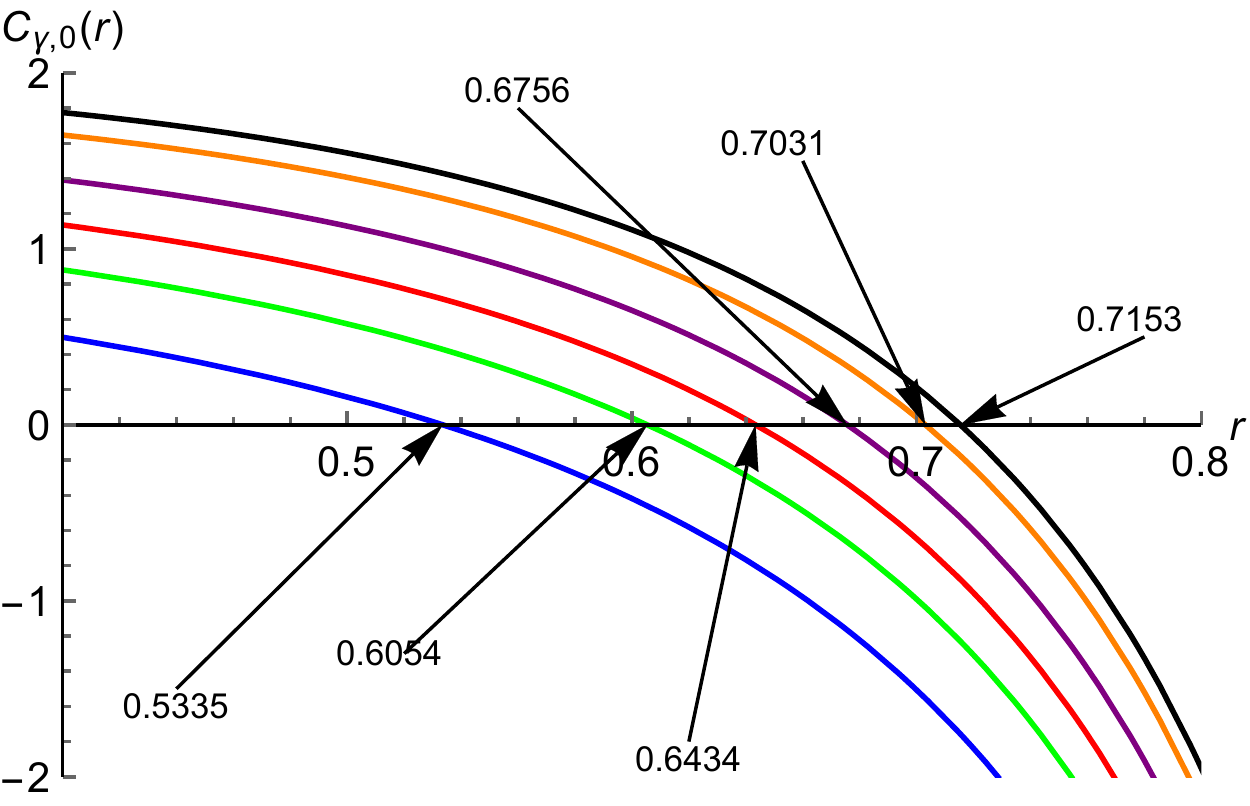}
		\,
		\includegraphics[width=0.48\linewidth]{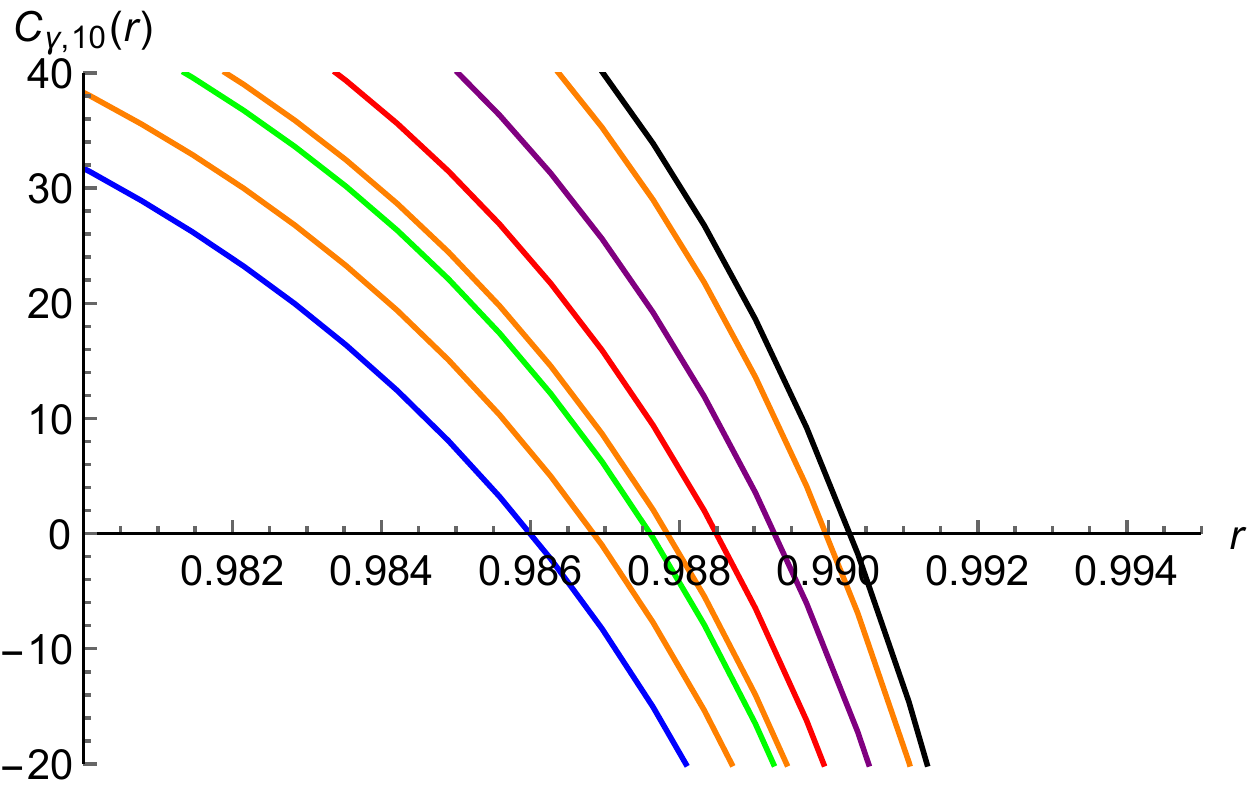}
	\end{center}
	\caption{The graph of $C_{\gamma,0}(r) $ and $C_{\gamma,10}( r)$ in $ (0,1) $ when $\gamma=0,0.3,0.5,0.7,0.9,1$.}
	\label{figure-4.2-d}
\end{figure}
\begin{figure}[!htb]
	\begin{center}
		\includegraphics[width=0.48\linewidth]{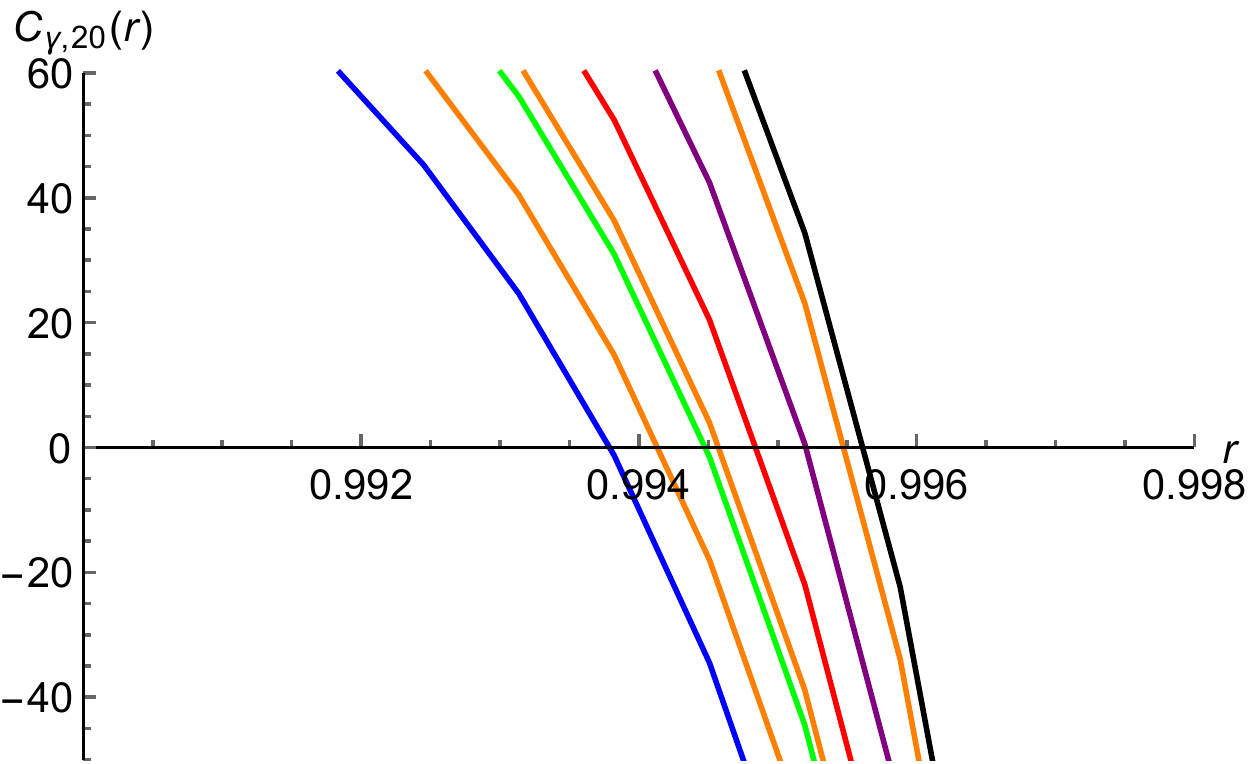}
		\,
		\includegraphics[width=0.48\linewidth]{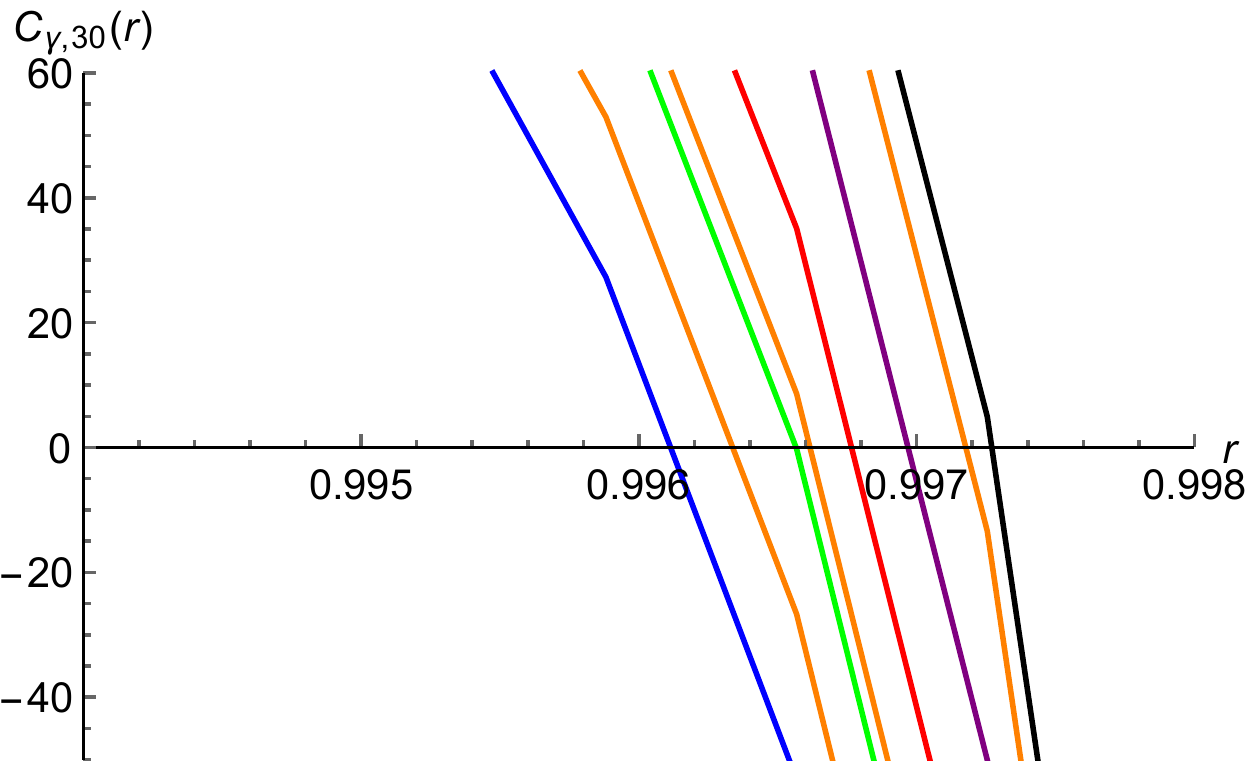}
	\end{center}
	\caption{The graph of $C_{\gamma,20}(r) $ and $C_{\gamma,30}( r)$ in $ (0,1) $ when $\gamma=0,0.3,0.5,0.7,0.9,1$.}
	\label{figure-4.2-e}
\end{figure}

From Table \ref{tabel-4.2-e}, for fixed values of $\alpha$, we observe that Bohr radius $R(\gamma,\alpha)$ is monotonic increasing in $\gamma \in [0,1)$. In the Table \ref{tabel-4.2-e}, the notation $(R(\gamma_{1},\alpha) \nearrow R(\gamma_{2},\alpha)]$ means that the value of $R(\gamma,\alpha)$ is monotonically increasing from $\lim _{\gamma \rightarrow \gamma ^{+}_{1}}= R(\gamma_{1},\alpha)$  to $ R(\gamma_{2},\alpha)$ when $\gamma_{1} < \gamma \leq \gamma_{2}$. The Figure \ref{figure-4.2-d} and Figure \ref{figure-4.2-e} are devoted to the graphs of $C_{\gamma,\alpha}(r)$ for various values of  $\gamma$ and $\alpha$. 
\begin{cor} \label{him-P7-cor-2.14}
	Let $f : \Omega_{\gamma} \rightarrow \mathbb{D}$ be an analytic function with $f(z)= \sum_{n=0}^{\infty} a_{n}z^n$ in $\mathbb{D}$. Then for $\alpha>-1$, the inequality \eqref{him-vasu-P6-e-2.70-a} holds
	for $|z|=r \leq R(\gamma,\alpha)$, where $R(\gamma,\alpha)$ is as in Theorem \ref{him-vasu-P6-thm-2.9}. In particular, for $\alpha=0$, we have
	\begin{equation} \label{him-vasu-P6-e-2.78}
	\mathcal{C}^{0}_{f}(r) \leq \frac{1}{r}\, \ln \left(\frac{1}{1-r}\right)
	\end{equation}
	for $|z|=r \leq R_{0}(\gamma)$, where $R_{0}(\gamma)$ is the smallest root in $(0,1)$ of $C_{\gamma}(r)=0$, where 
	\begin{equation*}
	C_{\gamma}(r)=(3+\gamma)(1-r)\ln\left(\frac{1}{1-r}\right) -2r.
	\end{equation*}
	The constant $R_{0}(\gamma)$ cannot be improved further. 
\end{cor}
For $\Omega_{\gamma}=\mathbb{D}$ {\it i. e.,} $\gamma=0$, using Corollary \ref{him-P7-cor-2.14}, we obtain the Bohr inequality for the Ces\'{a}ro operator for analytic functions $f:\mathbb{D} \rightarrow \mathbb{D}$.  
\begin{cor} \label{him-P7-cor-2.14-a}
	Let $f : \mathbb{D} \rightarrow \mathbb{D}$ be an analytic function with $f(z)= \sum_{n=0}^{\infty} a_{n}z^n$ in $\mathbb{D}$. Then for $\alpha>-1$, the inequality \eqref{him-vasu-P6-e-2.70-a} holds
	for $|z|=r \leq R(0,\alpha)$, where $R(0,\alpha)$ is as in Theorem \ref{him-vasu-P6-thm-2.9}. In particular, for $\alpha=0$, we have
	\begin{equation} \label{him-vasu-P6-e-2.78-a}
	\mathcal{C}^{0}_{f}(r) \leq \frac{1}{r}\, \ln \left(\frac{1}{1-r}\right)
	\end{equation}
	for $|z|=r \leq R_{0}(0)$, where $R_{0}(0)$ is the smallest root in $(0,1)$ of $C_{0}(r)=0$, where 
	\begin{equation*}
	C_{0}(r)=3(1-r)\ln\left(\frac{1}{1-r}\right) -2r.
	\end{equation*}
	The constant $R_{0}(0)$ cannot be improved further. 
\end{cor}
\vspace{3mm}

\section{Bohr inequality for Bernardi operator}
In similar fashion to the Bohr-type radius problem for the operator-valued $\alpha$-Ces\'{a}ro operator, we also study the Bohr-type radius problem for the operator-valued Bernardi operator. For $f:\mathbb{D} \rightarrow \mathcal{B}(\mathcal{H})$ analytic function with $f(z)= \sum_{n=m}^{\infty} A_{n}z^n$ in $\mathbb{D}$, where $A_{n} \in \mathcal{B}(\mathcal{H})$ for all $n \geq m$ and $m\geq 0$ is an integer with $\beta >-m$, we define Bernardi operator by 
\begin{equation*}
L_{\beta,\mathcal{H}}[f](z):= (1+\beta) \sum\limits_{n=m}^{\infty} \frac{A_{n}}{n+\beta} z^n=(1+\beta) \int\limits_{0}^{1}f(zt)\,t^{\beta-1}\,dt \,\,\,\,\,\, \mbox{for}\,\,\, z \in \mathbb{D}.
\end{equation*}
For $f(z)= \sum_{n=m}^{\infty} a_{n}z^n$ is complex-valued analytic function in $\mathbb{D}$, $L_{\beta,\mathcal{H}}$ reduces to complex-valued Bernardi operator $L_{\beta}$ (see \cite{Miller-Mocanu-diff-book}). For $\beta=1$ and $m=0$, we obtain the well-known Libera operator (see \cite{Miller-Mocanu-diff-book}) defined by 
\begin{equation*}
L[f](z):= 2\sum\limits_{n=0}^{\infty} \frac{a_{n}}{n+1} z^n=2\int\limits_{0}^{1}f(tz)\,dt \,\,\,\,\,\, \mbox{for}\,\,\, z \in \mathbb{D}.
\end{equation*}
For $\beta=0$, $m=1$, and $g(z)=\sum_{n=1}^{\infty} b_{n}z^n$, we obtain the well-known Alexander operator (see \cite{duren-1983}) defined by
\begin{equation*}
J[g](z):=\int\limits_{0}^{1} \frac{g(tz)}{t}\,dt =\sum\limits_{n=1}^{\infty} \frac{b_{n}}{n} z^n \,\,\,\,\,\, \mbox{for}\,\,\, z \in \mathbb{D},
\end{equation*} 
which has been extensively studied in the univalent function theory.
\par
In this section we study Bohr inequality for Barnardi operator $L_{\beta, \mathcal{H}}[f]$ when analytic functions $f :\Omega \rightarrow \mathcal{B}(\mathcal{H})$ for $f(z)=\sum_{n=m}^{\infty}A_{n}z^n$ in $\mathbb{D}$. 
Before going to establish Bohr inequality for the operator $L_{\beta, \mathcal{H}}$, we prove the following results, which are more general versions of Theorem \ref{him-P7-thm-1.3} and Theorem \ref{him-vasu-P6-thm-2.8}.
\begin{thm} \label{him-P7-thm-1.7}
	Fix $m \in \mathbb{N} \cup \{0\}$. Let $f \in H^{\infty}(\Omega,\mathcal{H})$ be given by $f(z)= \sum_{n=m}^{\infty} A_{n}z^n$ in $\mathbb{D}$ with $\norm{f(z)}_{H^{\infty}(\Omega,\mathcal{H})} \leq 1$, where $A_{m}=\alpha_{m}I$ for $|\alpha_{m}|<1$ and $A_{n} \in \mathcal{B}(\mathcal{H})$ for all $n \geq m$. If $\phi=\{\phi_{n}(r)\}_{n=m}^{\infty} \in \mathcal{G}$ satisfies the inequality
	\begin{equation} \label{him-P7-e-1.66}
	p \, \phi_{m}(r)> 2 \lambda_{\mathcal{H}} \sum_{n=m+1}^{\infty} \phi_{n}(r)\,\,\,\,\, \mbox{for} \,\,\,\, r \in [0, R_{\Omega}(p)),
	\end{equation}
	then the following inequality 
	\begin{equation} \label{him-P7-e-1.67}
	M_{f}(\phi, p,m, r):= \norm{A_{m}}^p \, \phi_{m}(r) + \sum_{n=m+1}^{\infty} \norm{A_{n}} \, \phi_{n}(r) \leq \phi_{m}(r)
	\end{equation}
	holds for $|z|=r \leq R_{\Omega}(p,m)$, where $R_{\Omega}(p,m)$ is the smallest root in $(0,1)$ of the equation 
	\begin{equation} \label{him-P7-e-1.68}
	p \, \phi_{m}(r)= 2 \lambda_{\mathcal{H}} \sum_{n=m+1}^{\infty} \phi_{n}(r).
	\end{equation}
\end{thm}
\begin{pf}
	Let $f \in H^{\infty}(\Omega,\mathcal{H})$ be of the form $f(z)= \sum_{n=m}^{\infty} A_{n}z^n$ in $\mathbb{D}$ with $\norm{f(z)}_{H^{\infty}(\Omega,\mathcal{H})} \leq 1$. Then we have $A_{m}=\alpha_{m}I$. We observe that $f(z)=z^m h(z)$, where $h:\Omega \rightarrow \mathcal{B}(\mathcal{H})$ is analytic function of the form $h(z)= \sum_{n=m}^{\infty} A_{n} z^{n-m}$ in $\mathbb{D}$ with $\norm{h(z)}_{H^{\infty}(\Omega,\mathcal{H})} \leq 1$. Then, in view of Definition \eqref{him-P7-e-1.24}, we have 
	\begin{equation} \label{him-P7-e-1.69}
	\norm{A_{n}} \leq \lambda_{\mathcal{H}}\norm{I- \abs{A_{m}}^2}= \lambda_{\mathcal{H}}\norm{I- \abs{\alpha_{m}}^2 I}=\lambda_{\mathcal{H}} (1- \abs{\alpha_{m}}^2)\, \,\,\,\, \mbox{for} \,\,\, n\geq m+1.
	\end{equation}
	Using \eqref{him-vasu-P6-e-2.59}, we obtain 
	\begin{align*}
	M_{f}(\psi, p,m, r) 
	& \leq \abs{\alpha_{m}}^p \, \psi_{m}(r) + \lambda_{\mathcal{H}} (1- \abs{\alpha_{m}}^2) \sum_{n=m+1}^{\infty} \psi_{n}(r) \\[2mm] 
	&= \psi_{m}(r) + \lambda_{\mathcal{H}} (1- \abs{\alpha_{m}}^2) \left(\sum_{n=m+1}^{\infty} \psi_{n}(r) - \frac{(1-|\alpha_{m}|^p)}{\lambda_{\mathcal{H}}(1-|\alpha_{m}|^2)} \psi_{m}(r)\right).
	\end{align*}
	Since $|\alpha_{m}|<1$, from the proof of Theorem \ref{him-P7-thm-1.3}, we have $(1-|\alpha_{m}|^p)/\lambda_{(\mathcal{H}}(1-|\alpha_{m}|^2)) \geq p/2\lambda_{\mathcal{H}}$ for $p \in (0,2]$, which leads to 
	$$
	M_{f}(\psi,p,m, r) \leq \psi_{m}(r) + \lambda_{\mathcal{H}} \left(1-\abs{\alpha_{m}}^2\right) \left(\sum_{n=1}^{\infty} \psi_{n}(r) - \frac{p}{2\lambda_{\mathcal{H}}} \psi_{m}(r)\right)
	$$
	and hence by \eqref{him-P7-e-1.66}, we obtain $M_{f}(\psi,p,m, r) \leq \psi_{m}(r)$ for $|z|=r \leq R_{\Omega}(p,m)$.
\end{pf}
\begin{thm} \label{him-P7-thm-1.8}
	Let $f \in H^{\infty}(\Omega_{\gamma},\mathcal{H})$ be of the form $f(z)= \sum_{n=m}^{\infty} A_{n}z^n$ in $\mathbb{D}$ with $\norm{f(z)}_{H^{\infty}(\Omega_{\gamma},\mathcal{H})} \leq 1$, where $A_{m}=\alpha_{m}I$ for $|\alpha_{m}|<1$ and $A_{n} \in \mathcal{B}(\mathcal{H})$ for all $n \geq m+1$. If $\phi=\{\phi_{n}(r)\}_{n=m}^{\infty} \in \mathcal{G}$ satisfies the following inequality
	\begin{equation} \label{him-P7-e-1.73}
	\phi_{m}(r)> \frac{2}{p(1+\gamma)} \sum_{n=m+1}^{\infty} \phi_{n}(r) \,\,\,\,\, \mbox{for} \,\,\,\, r\in [0,R(p,m,\gamma)),  
	\end{equation}
	then the inequality \eqref{him-P7-e-1.67} 
	holds for $|z|=r \leq R(p,m,\gamma)$, where $R(p,m,\gamma)$ is the smallest root in $(0,1)$ of the equation
	\begin{equation} \label{him-P7-e-1.74}
	\phi_{m}(x)=\frac{2}{p(1+\gamma)} \sum_{n=m+1}^{\infty} \phi_{n}(x).
	\end{equation}
	Moreover, when $\phi_{m}(x)<(2/p(1+\gamma)) \sum_{n=m+1}^{\infty} \phi_{n}(x)$ in some interval $\left(R(p,m,\gamma), R(p,m,\gamma)+ \epsilon\right)$ for $\epsilon >0$, then the constant $R(p,m,\gamma)$ cannot be improved further.
\end{thm}
\begin{pf}
	For $\Omega=\Omega_{\gamma}$, $\lambda_{\mathcal{H}}=1/(1+\gamma)$, the condition \eqref{him-P7-e-1.66} becomes
	\begin{equation*} 
	\phi_{m}(r)> \frac{2}{p(1+\gamma)} \sum_{n=m+1}^{\infty} \phi_{n}(r) \,\,\,\,\, \mbox{for} \,\,\,\, r\in [0,R(p,m,\gamma)),  
	\end{equation*}
	where $R(p,m,\gamma)$ is the smallest root in $(0,1)$ of the equation
	\begin{equation*}
	\phi_{m}(x)=\frac{2}{p(1+\gamma)} \sum_{n=m+1}^{\infty} \phi_{n}(x).
	\end{equation*}
	Then, by the virtue of Theorem \ref{him-P7-thm-1.3}, the required inequality \eqref{him-P7-e-1.67} holds for $r \in [0,R(p,m,\gamma))$. Our aim is to show that the radius $R(p,m,\gamma)$ cannot be improved further. That is, $\norm{A_{m}}^p \, \phi_{m}(r) + \sum_{n=m+1}^{\infty} \norm{A_{n}} \, \phi_{n}(r) > \phi_{m}(r)$ holds for any $r>R(p,m,\gamma)$ {\it i.e.,} for any $r \in \left(R(p,\gamma), R(p,\gamma)+ \epsilon\right)$. To show this, we consider the function $F_{a,m}: \Omega_{\gamma} \rightarrow \mathcal{B}(\mathcal{H})$ defined by $F_{a,m}(z)=z^m F_{a}(z)$, where $F_{a}$ is defined by \eqref{him-P7-e-1.36-a}. From the proof of Theorem \ref{him-P7-thm-1.3}, $\norm{F_{a}(z)}\leq 1$, and hence $\norm{F_{a,m}(z)}\leq 1$. Since $F_a(z)=A_0-\sum_{n=1}^{\infty}A_nz^n$ in $\mathbb{D}$, where $A_{0}, A_{n}$ are as in \eqref{him-P7-e-1.37-a}, then 
	\begin{equation*}
	F_{a,m}(z)=\left(\frac{a-\gamma}{1-a\gamma}\right)I z^m - (1-a^2) \sum_{n=m+1}^{\infty} \left(\frac{a^{n-m-1}(1-\gamma)^{n-m}}{(1-a\gamma)^{n-m+1}} \right)Iz^n\;\; \mbox{for}\;\; z\in\mathbb{D}.
     \end{equation*}	
	For the function $F_{a,m}$, we have 
	\begin{align} \label{him-P7-e-1.77} &
	\norm{A_{m}}^p \phi_{m}(r) + \sum_{n=m+1}^{\infty} \norm{A_{n}} \phi_{n}(r) \\[2mm] \nonumber
	&= \left(\frac{a-\gamma}{1-a\gamma}\right)^p \phi_{m}(r) + (1-a^2) \sum_{n=m+1}^{\infty} \frac{a^{n-1} (1-\gamma)^n}{(1-a\gamma)^{n+1}} \phi_{n}(r) \\[2mm] \nonumber
	&= \phi_{m}(r) + (1-a) \left(2\sum_{n=m+1}^{\infty} \phi_{n}(r) - p(1+\gamma) \phi_{m}(r)\right)+ \\[2mm] \nonumber
	& \,\,\,\,\, (1-a) \left(\sum_{n=m+1}^{\infty}\frac{a^{n-m-1} (1+a) (1-\gamma)^{n-m}}{(1-a\gamma)^{n-m+1}} \phi_{n}(r) - 2 \sum_{n=m+1}^{\infty} \phi_{n}(r)\right) + \\[2mm] \nonumber
	& \,\, \,\,\left(p(1+\gamma)(1-a) + \left(\frac{a-\gamma}{1-a\gamma}\right)^p -1\right) \phi_{m}(r) \\[2mm] \nonumber
	&= \phi_{m}(r) + (1-a) \left(2\sum_{n=m+1}^{\infty} \phi_{n}(r) - p(1+\gamma) \phi_{m}(r)\right) + O((1-a)^2)
	\end{align}
	as $a \rightarrow 1^{-}$. Also, we have that $$2\sum_{n=m+1}^{\infty} \phi_{n}(r) > p(1+\gamma) \phi_{m}(r)$$ for $r \in (R(p,m,\gamma), R(p,m,\gamma)+\epsilon)$. It is easy to see that the last expression of \eqref{him-P7-e-1.77} is strictly greater than $\phi_{m}(r)$ when $a$ is very close to $1$ {\it i.e.,} $a \rightarrow 1^{-}$ and $r \in (R(p,m,\gamma), R(p,m,\gamma)+\epsilon)$. This shows that the constant $R(p,m,\gamma)$ cannot be improved further. This completes the proof.
\end{pf}

Now we are in a position to establish Bohr inequality for Barnardi operator $L_{\beta, \mathcal{H}}[f]$ for analytic functions $f :\Omega_{\gamma} \rightarrow \mathcal{B}(\mathcal{H})$ of the form $f(z)=\sum_{n=m}^{\infty}A_{n}z^n$ in $\mathbb{D}$.
\begin{thm} \label{him-P7-thm-3.3}
Let $\beta >-m$ and $f : \Omega_{\gamma} \rightarrow \mathcal{B}(\mathcal{H})$ be an analytic function with $\norm{f(z)} \leq 1$ in $\Omega_{\gamma}$ such that $f(z)= \sum_{n=m}^{\infty} A_{n}z^n$ in $\mathbb{D}$, where $A_{m}=\alpha_{m}I$ for $|\alpha_{0}|<1$ and $A_{n} \in \mathcal{B}(\mathcal{H})$ for all $n \in m$, then
\begin{equation} \label{him-P7-e-3.8}
M_{\beta,\mathcal{H}}(r):=\sum_{n=m}^{\infty} \frac{\norm{A_{n}}}{n+\beta}\, r^n \leq \frac{1}{m+\beta} r^m
\end{equation}
for $|z|=r \leq R(m,\beta,\gamma)$, where $R(m,\beta,\gamma)$ is the smallest root in $(0,1)$ of $B_{m,\beta,\gamma}(r)=0$, where
\begin{equation} \label{him-P7-e-3.9}
B_{m,\beta,\gamma}(r)=\frac{2}{1+\gamma} \sum_{n=m+1}^{\infty} \frac{r^n}{n+\beta} - \frac{r^m}{m+\beta}.
\end{equation}
The constant $R(m,\beta,\gamma)$ is the best possible. 
\end{thm}
\begin{pf}
We note that $M_{\beta,\mathcal{H}}(r)$ can be expressed in the following form
\begin{equation*}
M_{\beta,\mathcal{H}}(r):=\sum_{n=m}^{\infty} \frac{\norm{A_{n}}}{n+\beta}\, r^n= \sum_{n=m}^{\infty} \norm{A_{n}} \, \phi_{n}(r)\,\,\,\,\,\,\, \mbox{with} \,\,\,\,\, \phi_{n}(r)=\frac{r^n}{n+\beta}
\end{equation*} 
and hence the condition \eqref{him-P7-e-1.73} becomes 
\begin{equation*}
\frac{r^m}{m+\beta} > \frac{2}{1+\gamma} \sum_{n=m+1}^{\infty} \frac{r^n}{n+\beta} \,\,\,\,\, \mbox{for} \,\,\,\, r\in [0,R(m,\beta,\gamma)), 
\end{equation*}
where $R(m,\beta,\gamma)$ is the smallest root of the equation \eqref{him-P7-e-3.9}. Now the inequality \eqref{him-P7-e-3.8} follows from Theorem \ref{him-P7-thm-1.8}. The sharpness of the constant $R(m,\beta,\gamma)$ follows from Theorem \ref{him-P7-thm-1.8}. This completes the proof.
\end{pf}
\begin{rem}
We observe that \eqref{him-P7-e-3.9} can also be written in the following form
\begin{equation*}
B_{m,\beta,\gamma}(r)=\frac{2}{1+\gamma} \sum_{n=1}^{\infty} \frac{r^{n+m}}{n+m+\beta} - \frac{r^m}{m+\beta}.
\end{equation*}
Thus, the root $R(m,\beta,\gamma)$ of $B_{m,\beta,\gamma}(r)=0$ is same as that of $L_{m,\beta,\gamma}(r)=0$, where 
\begin{equation} \label{him-P7-e-3.10}
L_{m,\beta,\gamma}(r)=\frac{2}{1+\gamma} \sum_{n=1}^{\infty} \frac{r^{n}}{n+m+\beta} - \frac{1}{m+\beta}.
\end{equation}
Therefore, \eqref{him-P7-e-3.10} yields that the roots of $L_{m,\beta,\gamma}(r)=0$ are same when the corresponding sums $m+\beta$ of $m$ and $\beta$ are the same. That is, for each fixed $i \in \mathbb{N}$, if $R(m_{i},\beta_{i},\gamma)$ is the root of $L_{m_{i},\beta_{i},\gamma}(r)=0$, then $R(m_{i},\beta_{i},\gamma)=R(m_{j},\beta_{j},\gamma)$ when $m_{i}+\beta_{i}=m_{j}+\beta_{j}$. For instance, $R(0,1,\gamma)=R(1,0,\gamma)=R(2,-1,\gamma)$, $R(0,2,\gamma)=R(1,1,\gamma)$.
\end{rem}
\begin{table}[ht]
	\centering
	\begin{tabular}{|l|l|l|l|l|}
		\hline
		$\gamma$& $R(0,1,\gamma)$& $R(0,2,\gamma)$& $R(1,2,\gamma)$& $R(4,0,\gamma)$ \\
		\hline
		$[0,0.2)$ & $[0.5828 \nearrow 0.6419)$& $[0.4742 \nearrow 0.5789)$ & $[0.4317 \nearrow 0.4833)$ & $[0.4090 \nearrow 0.4587)$\\
		\hline
		$[0.2,0.4)$& $[0.6419 \nearrow 0.6912)$& $[0.5289 \nearrow 0.5759)$& $[0.4833 \nearrow 0.5282)$& $[0.4587 \nearrow 0.5021)$\\
		\hline
		$[0.4,0.6)$& $[0.6912 \nearrow 0.7324)$& $[0.5759 \nearrow 0.6168)$& $[0.5282 \nearrow 0.5675)$& $[0.5021 \nearrow 0.5403)$\\
		\hline
		$[0.6,0.8)$& $[0.7324 \nearrow 0.7672)$& $[0.6168 \nearrow 0.6525)$& $[0.5675 \nearrow 0.6023)$& $[0.5403 \nearrow 0.5743)$\\
		\hline
		$[0.8,1)$& $[0.7672 \nearrow 0.7968)$& $[0.6525 \nearrow 0.6838)$& $[0.6023 \nearrow 0.6331)$& $[0.5743 \nearrow 0.6045)$\\
		
		\hline
	\end{tabular}
	\vspace{3mm}
	\caption{Values of $R(0,1,\gamma)$, $R(0,2,\gamma)$, $R(1,2,\gamma)$, and $R(4,0,\gamma)$ for various values of $\gamma \in [0,1)$.}
	\label{tabel-4.2-f}
\end{table}
\begin{figure}[!htb]
	\begin{center}
		\includegraphics[width=0.48\linewidth]{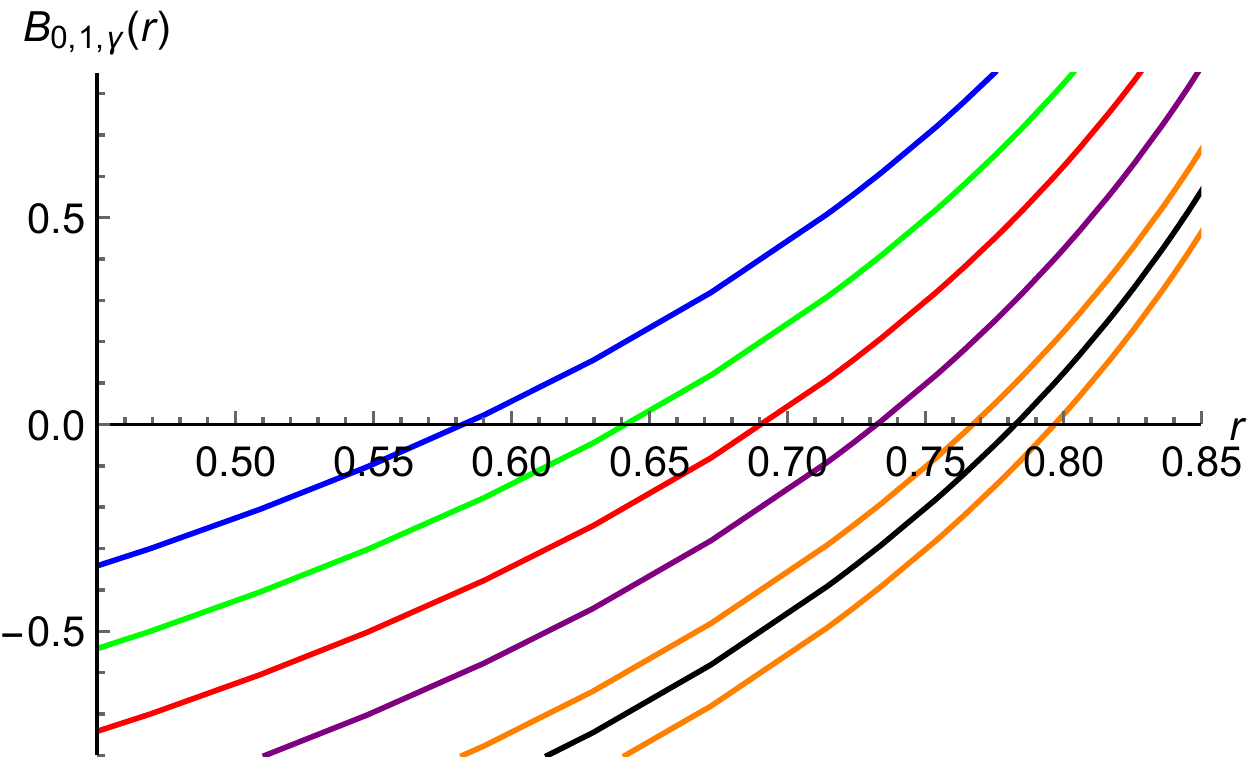}
		\,
		\includegraphics[width=0.48\linewidth]{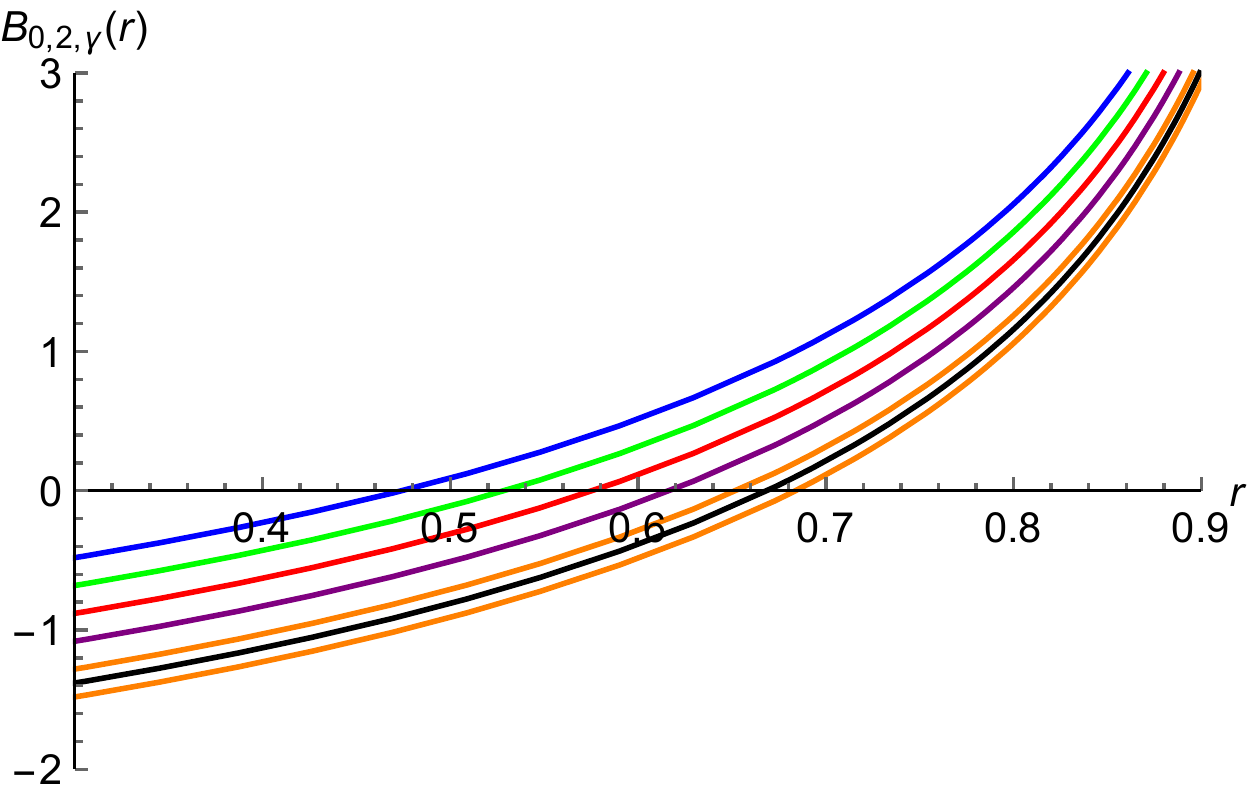}
	\end{center}
	\caption{The graph of $B_{0,1,\gamma}(r) $ and $B_{0,2,\gamma}( r)$ in $ (0,1) $ when $\gamma=0,0.2,0.4,0.6,0.8,0.9,1$.}
	\label{figure-4.2-f}
\end{figure}
\begin{figure}[!htb]
	\begin{center}
		\includegraphics[width=0.48\linewidth]{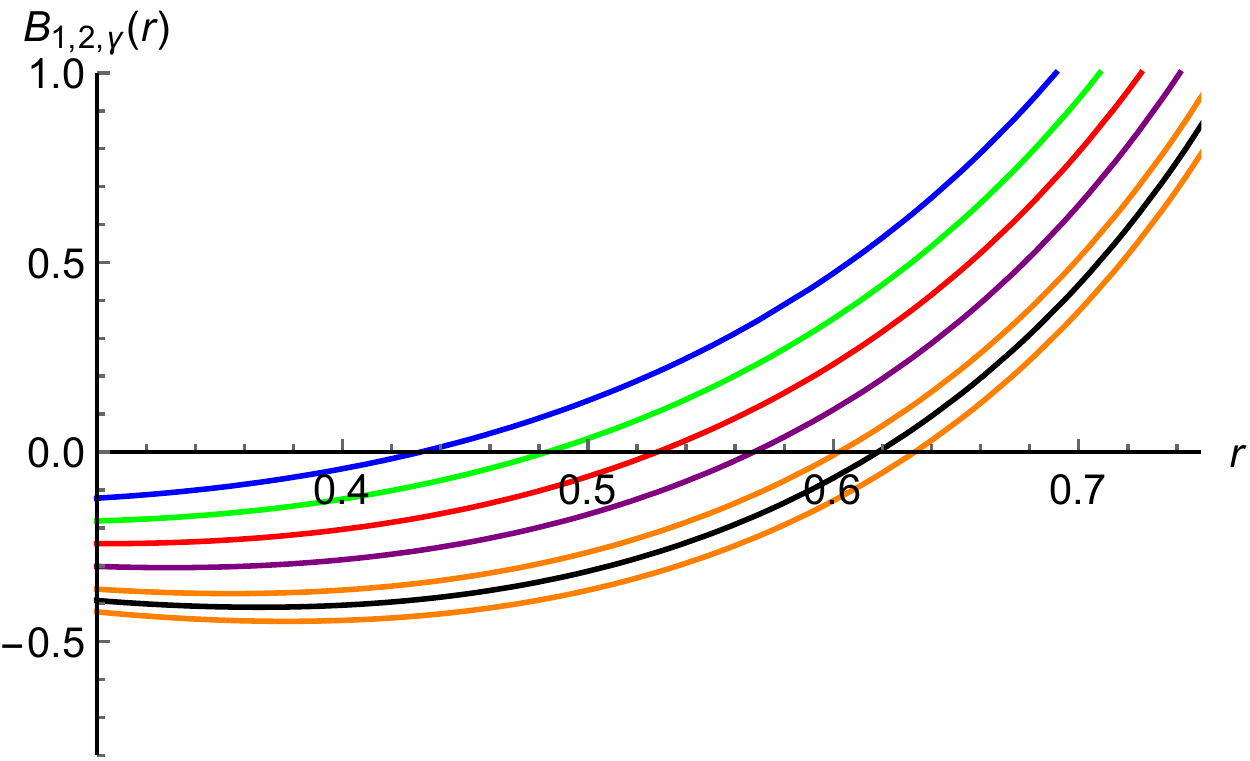}
		\,
		\includegraphics[width=0.48\linewidth]{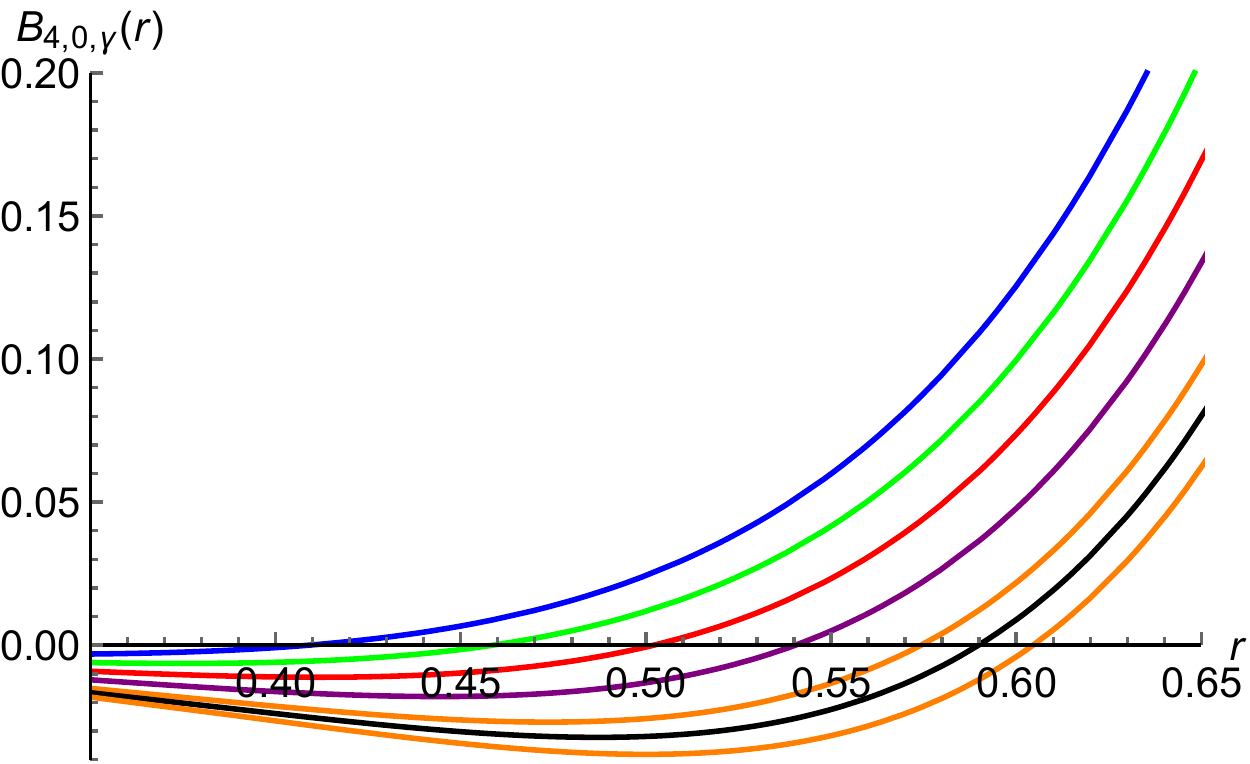}
	\end{center}
	\caption{The graph of $B_{1,2,\gamma}(r) $ and $B_{4,0,\gamma}( r)$ in $ (0,1) $ when $\gamma=0,0.2,0.4,0.6,0.8,0.9,1$.}
	\label{figure-4.2-g}
\end{figure}
From Table \ref{tabel-4.2-f}, for fixed values of $m$ and $\beta$, we observe that Bohr radius $R(m,\beta, \gamma)$ is monotonic increasing in $\gamma \in [0,1)$. In the Table \ref{tabel-4.2-f}, the notation $(R(m,\beta, \gamma_{1}) \nearrow R(m,\beta, \gamma_{2})]$ means that the value of $R(m,\beta, \gamma)$ is monotonically increasing from $\lim _{\gamma \rightarrow \gamma ^{+}_{1}} R(m,\beta, \gamma)= R(m,\beta, \gamma_{1})$  to $ R(m,\beta, \gamma_{2})$ when $\gamma_{1} < \gamma \leq \gamma_{2}$. The Figure \ref{figure-4.2-f} and Figure \ref{figure-4.2-g} are devoted to the graphs of $B_{m,\beta,\gamma}(r)$ for various values of  $m, \beta$ and $\gamma$.
\begin{cor}
Let $f$ be as in Theorem \ref{him-P7-thm-3.3} with $m=0$ and $\beta=1$. Then 
\begin{equation*}
\sum_{n=0}^{\infty} \frac{\norm{A_{n}}}{n+1}\, r^n \leq 1
\end{equation*}
for $|z|=r \leq R(0,1,\gamma)$, where $R(0,1,\gamma)$ is the smallest root in $(0,1)$ of 
\begin{equation} \label{him-P7-e-3.11}
\frac{2}{1+\gamma} \sum_{n=1}^{\infty} \frac{r^n}{n+1} = 1.
\end{equation}
The constant $R(0,1,\gamma)$ is the best possible. 
\end{cor}
\begin{cor}
	Let $f$ be as in Theorem \ref{him-P7-thm-3.3} with $m=1$ and $\beta=0$. Then 
	\begin{equation*} 
	\sum_{n=1}^{\infty} \frac{\norm{A_{n}}}{n}\, r^n \leq  r
	\end{equation*}
	for $|z|=r \leq R(1,0,\gamma)$, where $R(1,0,\gamma)$ is the smallest root in $(0,1)$ of 
	\begin{equation} \label{him-P7-e-3.13}
	\frac{2}{1+\gamma} \sum_{n=2}^{\infty} \frac{r^n}{n} =r.
	\end{equation}
	The constant $R(1,0,\gamma)$ is the best possible. 
\end{cor}

\noindent\textbf{Acknowledgment:} 
The first author is supported by SERB-CRG and the second author is supported by CSIR (File No: 09/1059(0020)/2018-EMR-I), New Delhi, India.

\end{document}